\documentclass[]{article}
\usepackage[utf8]{inputenc} 
\usepackage[T1]{fontenc}    
\usepackage{hyperref}       
\usepackage{url}            
\usepackage{microtype}
\usepackage{graphicx}
\usepackage{array}
\usepackage{color} 
\usepackage{tabularx}
\usepackage{amsmath,mathdots,amsthm}
\usepackage{amssymb}
\usepackage{amsfonts}
\usepackage{txfonts}
\usepackage{arxiv}
\usepackage{xcolor}
\usepackage{bm}
\usepackage{soul}
\usepackage{float}
\usepackage{cite}
\usepackage{booktabs, multirow}

\hypersetup{
	colorlinks=true,                          
	linkcolor=blue, 
	citecolor=blue, 
	urlcolor=black  
} 
\usepackage[colorinlistoftodos]{todonotes}
\usepackage[verbose]{placeins}

\newtheorem*{rmrk*}{Remark}

\newtheorem*{definition*}{Definition}

\newtheorem*{thm*}{Theorem}

\newtheorem*{prp*}{Proposition}

\title{An exactly curl-free finite-volume scheme for a hyperbolic compressible barotropic two-phase model}

\author{
	\href{https://orcid.org/0000-0002-7812-3161}{\includegraphics[scale=0.06]{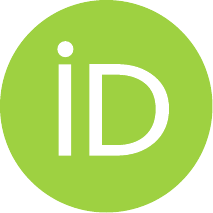}\hspace{1mm}Laura R\'io-Mart\'in}
	\thanks{corresponding author.\newline Email addresses: \texttt{laura.delrio@unitn.it} (L. R\'io-Mart\'in),
		\texttt{firas.dhaouadi@unitn.it} (F. Dhaouadi)
		\texttt{michael.dumbser@unitn.it} (M. Dumbser)} $^{\;a}$, 
	\hspace{3mm}
	\href{https://orcid.org/0000-0002-7150-3313}{\includegraphics[scale=0.06]{orcid.pdf}\hspace{1mm}Firas Dhaouadi$^b$}
	\hspace{3mm}
	\href{https://orcid.org/0000-0002-8201-8372}{\includegraphics[scale=0.06]{orcid.pdf}\hspace{1mm}Michael Dumbser$^b$}
}

\affiliation{	$^a$ Department of Information Engineering and Computer Science, University of Trento, via Sommarive 9, Povo, 38123 Trento, Italy.
\\
$^b$ Department of Civil, Environmental and Mechanical Engineering, University of Trento, Via Mesiano 77, 38123 Trento, Italy.}
\date{}

\renewcommand{\headeright}{L. R\'io-Mart\'in, F. Dhaouadi AND M. Dumbser}
\renewcommand{\undertitle}{}
\renewcommand{\shorttitle}{}

\hypersetup{
	pdftitle={},
	pdfsubject={},
	pdfauthor={Laura del Rio Martin, Firas Dhaouadi, Michael Dumbser},
	pdfkeywords={First keyword, Second keyword, More},
}

\newcommand{\pa}{^{I}}
\newcommand{\pb}{^{I\!I}}
\renewcommand{\u}[1]{u_{#1}} 
\newcommand{\uI}[1]{u\pa_{#1}} 
\newcommand{\uII}[1]{u\pb_{#1}} 
\newcommand{\w}[1]{w_{#1}} 
\newcommand{\bu}{\boldsymbol{u}} 
\newcommand{\buI}{\boldsymbol{u}\pa} 
\newcommand{\buII}{\boldsymbol{u}\pb} 
\newcommand{\bw}{\boldsymbol{w}} 
\newcommand{\bx}{\boldsymbol{x}} 
\newcommand{\bn}{\boldsymbol{n}} 
\newcommand{\pI}{p\pa}
\newcommand{\pII}{p\pb}
\newcommand{\cI}{c\pa}
\newcommand{\cII}{c\pb}
\newcommand{\aI}{a\pa}
\newcommand{\aII}{a\pb}
\renewcommand{\beta}{\boldsymbol{\eta}} 

\newcommand{\alphaI}{\alpha\pa}
\newcommand{\alphaII}{\alpha\pb}
\newcommand{\muI}{\mu\pa}
\newcommand{\muII}{\mu\pb}
\newcommand{\rhoI}{\rho\pa}
\newcommand{\rhoII}{\rho\pb}

\newcommand{\bg}{\boldsymbol{g}} 
\newcommand{\eI}{e\pa} 
\newcommand{\eII}{e\pb} 
\newcommand{\gammaI}{\gamma\pa}
\newcommand{\gammaII}{\gamma\pb}

\newcommand{\bQ}{\mathbf{Q}}

\newcommand{\bF}{\mathbf{F}}

\newcommand{\bB}{\mathbf{B}}
\newcommand{\bS}{\mathbf{S}}

\newcommand{\bFb}{\mathbf{F}^b}

\newcommand{\bBb}{\mathbf{B}^b}
\newcommand{\bSb}{\mathbf{S}^b}
\newcommand{\bFv}{\mathbf{F}^v}

\newcommand{\bGv}{\mathbf{G}^v}
\newcommand{\bBv}{\mathbf{B}^v}

\newcommand{\bDb}{\mathbf{D}^b}

\newcommand{\x}[1]{x_{#1}}

\newcommand{\ddx}[1]{\mathrm{d}{#1}}

\newcommand{\y}[1]{y_{#1}}

\newcommand{\pd}[2]{\frac{\partial #1}{\partial #2}}
\newcommand{\pdt}[1]{\frac{\partial #1}{\partial t}}

\newcommand{\prn}[1]{\left( #1 \right)}
\newcommand{\calE}{\mathcal{E}}

\begin{document}
	\maketitle

	\begin{abstract}
	We present a new second order accurate structure-preserving finite volume scheme for the solution of the compressible barotropic two-phase model of Romenski \textit{et. al} \cite{Romenski2007TwoPhase,Romenski2010TwoPhase} in multiple space dimensions. The governing equations fall into the wider class of symmetric hyperbolic and thermodynamically compatible (SHTC) systems and consist of a set of first-order hyperbolic partial differential equations (PDE). In the absence of algebraic source terms, the model is subject to a curl-free constraint for the relative velocity between the two phases. The main objective of this paper is therefore to preserve this structural property exactly also at the discrete level. The new numerical method is based on a staggered grid arrangement where the relative velocity field is stored in the cell vertexes, while all the remaining variables are stored in the cell centers. This allows the definition of discretely compatible gradient and curl operators which ensure that the discrete curl errors of the relative velocity field remain zero up to machine precision. A set of numerical results confirms this property also experimentally.
	\end{abstract}

\keywords{Hyperbolic equations \and compressible multiphase flows \and curl-free schemes \and finite-volume schemes \and staggered grids}

\section{Introduction}
Multiphase flows are ubiquitous in nature and in our everyday life. Already the rather simple flow of water with moving free surface actually involves both the dynamics of the moving liquid and of the surrounding air. The modeling and numerical simulation of such multi-material systems remains challenging even nowadays, as no prevalent approach seems to be efficient for all applications.

In this paper, we are concerned with compressible barotropic two-phase flows and in particular we consider the model of Romenski forwarded in \cite{Romenski1998,Romenski2004,Romenski2007TwoPhase,Romenski2010TwoPhase}. The model belongs to the so-called class of Symmetric Hyperbolic and Thermodynamically Compatible (SHTC) systems discovered by Godunov and Romenski in a series of seminal works \cite{Godunov1961,Romenski1998,Godunov2003} and which was extensively studied and solved numerically more recently in \cite{Romenski2004,Romenski2007TwoPhase,Romenski2010TwoPhase,Peshkov2018SHTC, Thein2022,Lukacova2023,Rio-Martin2023,Thomann2023,Thomann2024}. Some of the reasons motivating our interest in this particular model are the fact that it is first order symmetric hyperbolic, it can be rewritten in a similar form of the well-known Baer-Nunziato model \cite{BaerNunziato1986} and offers a rather general framework that can be extended also to viscous and heat conducting multi-material flows, including mixtures of solids, liquids and gases within one and the same mathematical framework. 
One particular property of the model is that under certain assumptions, one of its equations (describing the evolution of the relative velocity) is bound by a curl-free constraint. The latter comes out as an involution in the sense that it is rather a direct consequence of the main system of equations and not an additional condition to be supplied. The most prominent examples of involutions are the well-known divergence-free condition of the magnetic field in the Maxwell and MHD equations, or the curl-free property of the deformation gradient in solid mechanics. Such stationary differential constraints (involutions) are present in many other systems of physics and continuum mechanics, see for example \cite{GodunovRomenski1972,Godunov1972MHD,Godunov2003,Alic2009,Brown2012,Peshkov2018SHTC,Dumbser2020GLM,Busto2021HyperbolicDispersion,Chiocchetti2021,Dhaouadi2022NSK,Dhaouadi2023Heat}. Involution constraints generally are of little consequence at the continuous level, since solutions of the governing PDE system obey them by definition. One cannot say the same at the discrete level, as extra attention needs to be paid for discrete solutions to remain compatible with the involution, thus preserving their physical relevancy. This spurred the development of particular numerical methods allowing to preserve differential involutions, mainly divergence and curl constraints, at the discrete level. Examples of such methods include for example constrained transport methods \cite{Yee1966,Brecht1981,Holland1983,Evans1988,Devore1991,Dai1998,Toth2000,Gardiner2005,Boscheri2023MHD}, divergence and curl cleaning approaches  \cite{Munz2000,Dedner2002,Dedner2003,Jacobs2009,Dumbser2020GLM,Chiocchetti2021,Busto2021HyperbolicDispersion,Dhaouadi2022NSK} and structure-preserving discretizations \cite{Balsara1999,Balsara2004,Torrilhon2004,JeltschTorrilhon2006,Xu2016,Hazra2019,Boscheri2021SIGPR,Balsara2023CurlFree,Boscheri2023SPDG}.

In this context, we present the following contribution, the main objective being to solve numerically the system of two-phase flows forwarded by Romenski while guaranteeing that its intrinsic curl-free constraint is respected exactly at the discrete level. For that, we will consider the second-order structure-preserving finite volume scheme proposed in~\cite{Boscheri2021SIGPR} and further studied in~\cite{Dumbser2020Curl,Peshkov2021NN,Dhaouadi2023NSK,Chiocchetti2023}. One of the novelties here is that the constrained vector field is a velocity field, thus requiring extra effort in the implementation of compatible boundary conditions as soon as reflective walls are considered. 

This paper is organized as follows. In section~\ref{sec:model}, we provide the notations and recall the structure of the compressible two-phase model under consideration. We present the equations of state that shall be used in the simulations and offer a concise review of the hyperbolic nature of the equations. In section \ref{sec:Numerical_method}, we explain all parts of the proposed curl-free numerical scheme. We show how the discretization of the system is split between a primary and a dual grid, and we recall how this staggering allows us to recover the curl-free property exactly at the discrete level. In section \ref{sec:Numerical_Results} we present some numerical results. We show that the proposed scheme exhibits second-order convergence for a smooth vortex-type solution. We  compare the approximate solution calculated with the proposed methodology and a reference solution for Riemann problems in one and two space dimensions. In particular, we show comparisons with reference solutions for a one-dimensional Riemann problem and for a radial explosion test. The computed solution for a dam-break problem will also be compared with the reference solution calculated using an equivalent Baer-Nunziato model. In this test case, we will describe how the compatible wall boundary conditions have been implemented. Finally, the performance of the scheme to satisfy the curl-free constraint of the relative velocity is illustrated via the simulation of a Kelvin-Helmholtz instability, which exhibits rather complex flow features. In all cases, plots of the curl errors over time are given as evidence of the structure-preserving property of the scheme.     

\section{The two-phase model of Romenski \textit{et al.}}
\label{sec:model}
\subsection{Governing equations}
We consider the model of Romenski proposed in \cite{Romenski1998,Romenski2001,Romenski2004,Romenski2007TwoPhase,Romenski2010TwoPhase}, describing the motion of a multiphase medium formed by the mixture of two compressible fluids. Here, we neglect any effects due to viscosity and inter-phase friction,  and we assume that the motion takes place in the absence of any heat or mass exchange. As to clarify the notations in what follows, a superscript shall be used to designate the phase ($I$ for phase one, $I\! I$ for phase two). Subscripts will be reserved for vector components and matrix entries. Repeated subscript summation via the usual Einstein summation convention is implied. Under these assumptions and notations, the equations of motions are given as follows:
\begin{subequations}
	\begin{align}
		&\frac{\partial \alphaI}{\partial t}+u_k \frac{\partial \alphaI}{\partial x_{k}}=0   ,\label{eq:cons_volume_frac}\\
		&\frac{\partial\alphaI\rhoI}{\partial t}+\frac{\partial(\alphaI\rhoI 	\uI{k})}{\partial x_{k}}=0,\label{eq:cons_mass_density1}\\
		&\frac{\partial\alphaII\rhoII}{\partial t}+\frac{\partial(\alphaII\rhoII 	\uII{k})}{\partial x_{k}}=0, \qquad \alphaII = 1-\alphaI,\label{eq:cons_mass_density2}\\
		&\frac{\partial\rho \u{i}}{\partial t}+ \frac{\partial\prn{\rho \u{i} u_k+\Pi_{ik}}}{\partial x_{k}}=g\textcolor{red}{_i} \rho, \quad \Pi_{ik} = p\delta_{ik} + \rho w_i \, E_{w_k} , \quad i\in\{1,\ldots d\}, \label{eq:momentum_ini}\\
		&\frac{\partial w_k}{\partial t}+\frac{\partial (w_l u_l + \varphi)}{\partial x_{k}}+ u_l\left(\frac{\partial w_k}{\partial x_{l}}-\frac{\partial w_l}{\partial x_{k}} \right)=0, \quad k\in\{1,\ldots d\},
		\label{eq:relative_vel_ini}
	\end{align}
	\label{eq:TwoFluid}
\end{subequations}  
Here, $t\in\mathbb{R}_+$ is the time and $\mathbf{x}\in\mathbb{R}^d$ is the space coordinates vector. The quantities $\alpha^j$, $\rho^j$ and $\bu^j=(u^j_1,\cdots,u^j_d)$, where $j\in\{I,I\!I\}$, are the phase average volume fraction, density and velocity field, respectively, of the $j^{th}$ component. The mixture density $\rho$, the mixture velocity $\bu$, and the relative velocity $\mathbf{w}$ are then given by
\begin{equation*}
	\rho = \alphaI \rhoI  + \alphaII \rhoII, \quad \bu = \frac{\alphaI \rhoI \buI + \alphaII \rhoII \buII}{\rho}, \quad \bw = \buI - \buII.
\end{equation*}	
Note that the system's state is fully determined by the knowledge of the variables $\left\{ \alphaI, \rhoI, \rhoII, \bu, \bw \right\}$ and any other quantity should be understood as a function of the latter and not as an independent degree of freedom. 
Under the previous definitions, equation \eqref{eq:cons_volume_frac} describes the transport of the volume fraction. Equations (\ref{eq:cons_mass_density1},\ref{eq:cons_mass_density2}) describe mass conservation for each phase. Equation \eqref{eq:momentum_ini} is the mixture momentum conservation equation. Note that the vector $\bg$ is the gravity field, which we will only consider in some test cases. Lastly, equation \eqref{eq:relative_vel_ini} is the balance law for the relative velocity. Note that in the absence of source terms in the latter, the equation is subject to a curl-free constraint, provided the field is initially as such, that is
\begin{equation}
	\text{if}	\quad \nabla \times \bw = 0, \quad \text{at} \ t=0 \quad \text{then} \quad \nabla \times \bw = 0 \quad \forall t \geq 0,
	\label{eq:curl=0}
\end{equation}  
and which is an immediate consequence of equation \eqref{eq:relative_vel_ini}. For this reason, the field $\bw$ can be seen as the gradient of a scalar function $\phi$. The total energy density of this system is defined as
\begin{equation*}
	\calE(\alphaI,\rhoI,\rhoII,\bw,\bu) = \frac{1}{2} \rho u_{l}u_{l} + \rho E(\alphaI,\rhoI,\rhoII,\bw),
\end{equation*}
where $E(\alphaI,\rhoI,\rhoII,\mathbf{w})$ is the specific internal energy of the mixture and which can be written in separable form as  
\begin{equation*}
	E(\alphaI,\rhoI,\rhoII,\mathbf{w}) = \cI \eI(\rhoI) + \cII \eII(\rhoII) + \frac{\cI \cII}{2} \w{l}\w{l},
\end{equation*}
where $e^j$ is the internal specific energy of the $j^{th}$ phase, whose expression is given by the equation of state of the corresponding phase. The mass fractions $c^j$ are only introduced to ease notation and are given as functions of the energy functional arguments as follows 
\begin{equation*}
	\cI(\rhoI,\rhoII,\alphaI) = \frac{\alphaI \rhoI}{\alphaI \rhoI  + \prn{1-\alphaI} \rhoII }, \quad 	\cII(\rhoI,\rhoII,\alphaI) = \frac{\prn{1-\alphaI} \rhoII}{\alphaI \rhoI  + \prn{1-\alphaI} \rhoII}.
\end{equation*}
The tensor $\Pi_{ik}$ that appears in the mixture momentum equation \eqref{eq:momentum_ini} is the total mixture stress tensor, where the mixture pressure $p$ can be expressed as the average of the phase pressures so that 
\begin{equation*}
	p = \alphaI \pI + \alphaII \pII, \quad \text{where} \quad p^j = \prn{\rho^j}^2 \pd{e^j}{\rho^j}.
\end{equation*}
The scalar field $\varphi$ which appears in the relative velocity equation \eqref{eq:relative_vel_ini} writes as 
\begin{equation*}
	\varphi = \frac{\rho}{\alphaI}\pd{E}{\rhoI} - \frac{\rho}{\alphaII}\pd{E}{\rhoII} = \muI - \muII - \frac{\cI-\cII}{2} \w{l}\w{l}
\end{equation*}
where $\mu^j = e^j(\rho^j) + p^j/\rho^j$ is the chemical potential of the $j^{th}$ phase.
Finally, one can obtain an additional conservation law for the total energy of the mixture as a consequence of the system of equations \eqref{eq:TwoFluid}, which is written as
\begin{equation*}
	\pdt{\calE} + \pd{\prn{\calE u_{k} + \Pi_{ik} u_{k}  + \rho \varphi E_{\w{k}}  }}{x_k} = 0.
\end{equation*}

\subsection{Equations of state}	
Throughout this paper, we shall make use of either an ideal gas or a stiffened gas equation of state. In particular, the latter will be used whenever a liquid phase is considered. We recall hereafter the expressions of the internal energy and the corresponding pressure for a fluid of density $\rho$, in the barotropic case. For an ideal gas, we have   
\begin{equation}
	e(\rho)=\frac{\rho^{\gamma-1}}{ (\gamma-1)}, \quad p(\rho) = \rho^\gamma
	\label{eq:energy_ideal_gas}
\end{equation}
while for a stiffened gas, we have
\begin{equation}
	e(\rho)=\frac{c_0^2}{\gamma (\gamma-1)}\left(\frac{\rho}{\rho_0}\right)^{\gamma-1}+\frac{\rho_0 c_0^2 - \gamma p_0}{\gamma \rho}, \quad p(\rho)= p_0 + \frac{\rho_0 }{\gamma}c_0^2 \prn{\prn{\frac{\rho}{\rho_0}}^\gamma - 1}.
	\label{eq:energy_stiffened}
\end{equation}
In these expressions, $\gamma=c_p/c_V$ is the ratio of the heat capacity at constant pressure to heat capacity at constant volume, $p_0$, $\rho_0$ and $c_0$ are reference quantities for the pressure, density, and sound speed, respectively.

\subsection{Hyperbolicity}
The hyperbolicity of system~\eqref{eq:TwoFluid} was addressed in one dimension of space in~\cite{Romenski2004,Thein2022} while the multidimensional case was discussed for the first time in~\cite{Rio-Martin2023}. In particular, it was shown therein that in three dimensions of space, under considerations of rotational invariance, the system of equations~\eqref{eq:TwoFluid} has nine real eigenvalues in the $x-$ direction, whose expressions are recalled here
\begin{equation}
	\lambda_1 = \uII{1} - \aII, \quad \lambda_2 = \uI{1} - \aI, \quad \lambda_{3-7} = \u{1}, \quad 	\lambda_{8} = \uII{1} + \aII, \quad \lambda_{9} = \uI{1} + \aI,  
\end{equation}
where $a^j = \sqrt{\partial p^j/\partial \rho^j}$ is the phase averaged sound speed. In particular, the system of equations~\eqref{eq:TwoFluid} is only weakly hyperbolic for $d\geq2$ as the eigenvalue $\lambda=u_1$ is defective and lacks $d-1$ linearly independent eigenvectors. Such a shortcoming occurs in many systems of continuum mechanics that share a similar structure to the one considered here, namely systems of equations discussed in \cite{Godunov1961}, whenever a curl-free vector field is evolved in time, see for example \cite{Chiocchetti2021,Busto2021HyperbolicDispersion, Dhaouadi2022NSK,Dhaouadi2023Heat}. A known fix allowing to recover strong hyperbolicity in multiple dimensions of space is the addition of the well-known symmetrizing terms, also referred to as Godunov-Powell terms \cite{Godunov1972MHD,Powell1994MHD1,Powell1997MHD,Powell1999MHD}. In the case of the present system of equations, this would require modifying the momentum equation by adding the vector product $\prn{\nabla \times \bw}\times \rho E_{\bw}$ so that it becomes
\begin{equation*}
	\frac{\partial\rho \u{i}}{\partial t}+ \frac{\partial\prn{\rho \u{i} u_k+\Pi_{ik}}}{\partial x_{k}}  - \rho E_{w_k} \prn{\pd{\w{i}}{\x{k}} - \pd{\w{k}}{\x{i}}} = g^i \rho \label{eq:momentum_powell}.
\end{equation*}
Such a modification is only legitimate under the assumption that the initial data satisfies $\nabla \times \bw = 0$ as, in this case, the added term also cancels out in virtue of Equation~\eqref{eq:curl=0}. This actually restores the system's strong hyperbolicity, as shown in \cite{Rio-Martin2023}. A possible interpretation of this fact is that while the system of equations~\eqref{eq:TwoFluid} is weakly hyperbolic for general initial data, solutions of the associated initial-value problem are equivalent in some sense to the solutions of the strongly hyperbolic system, altered through the Godunov-Powell terms for a restricted set of well-prepared initial data satisfying $\nabla \times \bw = 0$ at $t=0$. An alternative treatment that allows to restore strong hyperbolicity without modifying the momentum conservation law was pointed out in \cite{Rio-Martin2023} and consists in the use of a generalized Lagrangian multiplier (GLM) technique, see \cite{Munz2000,Dedner2002,GassnerEntropyGLM,Busto2023new} for the original GLM method applied to divergence-type involutions and its more recent variant  \cite{Dumbser2020GLM,Busto2021HyperbolicDispersion,Chiocchetti2021,Dhaouadi2022NSK,Rio-Martin2023} for curl-type involutions. 

\section{Numerical method}
\label{sec:Numerical_method}
In this section, we propose a structure-preserving discretization to solve the system of equations~\eqref{eq:TwoFluid}. This discretization is based on a staggered mesh arrangement, where the relative velocity is stored in the vertices of the mesh, while the remaining variables, such as the mixture density, the momentum and the phase volume fractions are stored in the cell centers. The first part of this section is devoted to describing the used staggered grid, introducing also the necessary details of the notation. In addition, we will describe the compatible discrete operators for the gradient and the curl that are employed to discretize the relative velocity equation. Then, we introduce the flux splitting that is a key ingredient of our methodology. Finally, we provide a brief description of the second-order MUSCL-Hancock-type scheme used to discretize the remaining terms.

\subsection{Notation}
\begin{figure}
	\begin{center}
		\includegraphics[width=0.48\textwidth]{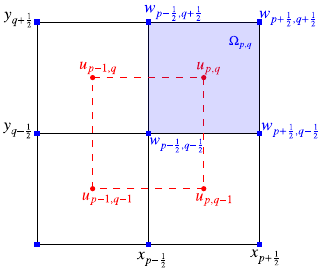} \qquad 
		\includegraphics[width=0.42\textwidth]{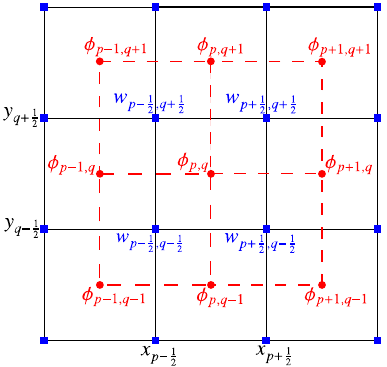}  
		\caption{Sketch of the staggered grid. The blue elements are the cells of the main grid $\Omega_{p,q}$, and those delimited by dashed red lines are the cells of the dual grid. Left: the relative velocity $\w{}$ is defined at the vertices of the main grid, and the rest of variables, such as $\u{}$, are defined in the cell centers of the control volumes. Right: The scalar field $\phi$ is also defined in the cell centers.}
		\label{fig:staggered_grid}
	\end{center}
\end{figure}

The main objective of this section is to present a compatible discretization method that satisfies the curl constraint for the relative velocity exactly also at the discrete level. With this purpose, we consider a discretization based on staggered grids, which will allow us to define some fields in the cell centers and others in the vertices of the elements of the main grid. Figure~\ref{fig:staggered_grid} shows the considered staggered grids. The blue elements are the cells of the main grid $\Omega_{p,q}$ and the elements limited by dashed red lines are the vertex-based dual cells $\Omega_{p\pm\frac{1}{2},q\pm\frac{1}{2}}$. As it is shown in the left sketch, the relative velocity $\w{}$ is defined in the vertices of the main grid, and the variables $\u{},\uI{},\uII{},\rho,\rhoI,\rhoII$ are defined in the cell centers. 
Below, we will provide a detailed description of the notation. Throughout this section, and to distinguish coordinate indices from time and space discretization indices, we will denote the former as $i,j,k$, the time discretization index as $n$ and the spatial discretization indices as $p,q$. For the sake of simplicity, we will describe the numerical method considering a two--dimensional domain $\Omega$, with $\x{1}=x$ and $\x{2}=y$. The computational domain $\Omega = [-\frac{L_{x}}{2},\frac{L_{x}}{2}]\times[-\frac{L_{y}}{2},\frac{L_{y}}{2}]$ is discretized with a uniform Cartesian grid composed of $N_{x}\times N_{y}$ cells. These cells are given by $\Omega_{p,q} = [x_{p-\frac{1}{2}}, x_{p+\frac{1}{2}}]\times [y_{q-\frac{1}{2}}, y_{q+\frac{1}{2}}]=[x_p-\frac{\Delta x}{2},x_p+\frac{\Delta x}{2}]\times[y_q-\frac{\Delta x}{2},y_q+\frac{\Delta y}{2}]$, with $(x_p,y_q)$ the discrete spatial coordinates located in the barycenter of the control volume $\Omega_{p,q}$, and $\Delta x=\frac{L_x}{N_x}$, $\Delta y=\frac{L_y}{N_y}$ the uniform mesh spacing in the $x-$ and $y-$directions.
The critical point of the presented structure-preserving method lies in the definition of a discrete gradient and curl, which will be defined at the staggered grid points and which will be compatible at the discrete level with the system structure. Below, the flux splitting considered for this problem is introduced, and the compatible operators are defined.

\subsection{Flux splitting}
Let $\bQ=\left(\alphaI, \alphaI\rhoI, \alphaII\rhoII, \rho \bu, \bw\right)^T$ be the state vector, then the PDE system~\eqref{eq:TwoFluid} can be written more compactly as
\begin{equation*}
	\partial_t \bQ + \nabla \cdot \bF(\bQ) + \bB(\bQ)\cdot\nabla\bQ =\bS(\bQ),
\end{equation*}
where $\bF(\bQ)$ is the nonlinear flux tensor, $\bB(\bQ)\cdot\nabla\bQ$ contains the non-conservative terms and $\bS(\bQ)$ is the algebraic source term. To apply the numerical method to solve the proposed PDE system, it is more appropriate to split the system into the form
\begin{equation}
	\partial_t \bQ + \nabla \cdot \left(\bFb(\bQ)+\bFv(\bQ)\right)+ \nabla \bGv(\bQ)+ \bBb(\bQ)\cdot\nabla\bQ + \bBv(\bQ)\cdot\nabla\bQ=\bSb(\bQ),
	\label{eq:flux_splitting}
\end{equation}
where the superscripts $b$ and $v$ distinguish terms that involve fields defined on the cell centers, or on the vertices, respectively. The terms appearing in  equation~\eqref{eq:flux_splitting} are defined as
\begin{align*}
	{\bFb} &=\begin{pmatrix}0\\  \alphaI\rhoI\u{k} \\ \alphaII\rhoII\u{k} \\\rho \u{i} \u{k}+p\delta_{ik} \\ 0\end{pmatrix}, \qquad 
	{\bFv}=\begin{pmatrix}0\\  \alphaI\rhoI\cII\w{k} \\ -\alphaII\rhoII\cI\w{k} \\\cI\cII\rho \w{i} \w{k} \\ 0\end{pmatrix}, \qquad 
	{\bGv}=\begin{pmatrix}0\\ 0 \\ 0 \\ 0\\\w{j} \u{j} + \varphi\end{pmatrix}, \\
	&{\bBb(\bQ)\cdot\nabla\bQ}=\begin{pmatrix}\u{k} \pd{\alphaI}{\x{k}}\\ 0 \\ 0 \\ 0 \\ 0\end{pmatrix}, \qquad
	{\bBv(\bQ)\cdot\nabla\bQ}=\begin{pmatrix}0\\ 0 \\ 0 \\ 0\\ \u{j}\left(\pd{\w{k}}{\x{j}}-\pd{\w{j}}{\x{k}} \right) \end{pmatrix}, \qquad	
	{\bSb}=\begin{pmatrix}0\\ 0 \\ 0 \\ g_i \rho \\ 0\end{pmatrix}, 
\end{align*}
where $\bFb$ is the flux containing the convective part as well as the pressure terms, $\bFv$ is the flux containing the terms related to the relative velocity, $\bGv$ comprises the components whose gradient will be calculated using the compatible discrete gradient operator that we will describe in the next section, $\bBb(\bQ)\cdot\nabla\bQ$ and $\bBv(\bQ)\cdot\nabla\bQ$ contain the non-conservative products related to the phase volume fraction and to the curl terms, respectively, and $\bSb$ are the source terms of the momentum equation related to gravity. 	
The subsystem 
\begin{equation}
	\partial_t \bQ + \nabla \cdot \left(\bFb(\bQ)+\bFv(\bQ)\right)+ \bBb(\bQ)\cdot\nabla\bQ =\bSb(\bQ),
	\label{eq:flux_splittingb}
\end{equation}
will be discretized explicitly using a second-order MUSCL-Hancock type finite volume method with a TVD limiter, using the Rusanov flux as approximate Riemann solver.
On the other hand, the discretization of 
\begin{equation}
	\partial_t \bQ + \nabla \bGv(\bQ)+ \bBv(\bQ)\cdot\nabla\bQ=\boldsymbol{0},
	\label{eq:flux_splittingv}
\end{equation}
will be performed using the compatible gradient and curl operators that will be described in the next section to obtain a curl-free method at the discrete level.

\subsection{Compatible gradient and curl operators}	

As we have indicated at the beginning of this section, the key point of the structure-preserving method presented in this paper lies in properly defining the discrete differential operators, considering the grid points we have just defined. Let $\phi_{p,q}^n=\phi(\x{p},\y{q},t^n)$ be a scalar field defined at the centers of the control volumes $\Omega_{p,q}$, the primitive variables ($\rho_{p,q}, \rhoI_{p,q},\rhoII_{p,q},\uI{p,q},\uII{p,q},\ldots$) which are stored in the cell centers ($\x{p}, \y{p}$) and the relative velocity field ($\w{p\pm\frac{1}{2},q\pm\frac{1}{2}} =\nabla\phi_{p\pm\frac{1}{2},q\pm\frac{1}{2}} $) that is stored at the dual grid points ($\x{p\pm\frac{1}{2}},\y{q\pm\frac{1}{2}}$). If there is no confusion, we will omit the superscript $n$ for the time level for the sake of simplicity. 
\\
The \textit{discrete gradient operator} $\nabla^h$ of the scalar field $\phi$ is defined in each vertex as a constant value that can be computed naturally using the finite differences of $\phi$ along the two independent directions via a central corner gradient (the right sketch of Figure~\ref{fig:staggered_grid} shows more details). The degrees of freedom of the compatible gradient operator read as
\begin{equation}
	\left(\nabla^h\phi^{h}\right)=\begin{pmatrix}
		\left(\partial_x^h\phi\right)_{p+\frac{1}{2},q+\frac{1}{2}}\\
		\left(\partial_y^h\phi\right)_{p+\frac{1}{2},q+\frac{1}{2}}
	\end{pmatrix}=
	\begin{pmatrix}
		\dfrac{1}{2}\dfrac{\phi_{p+1,q+1}-\phi_{p,q+1}+\phi_{p+1,q}-\phi_{p,q}}{\Delta x}
		\\[0.2cm]
		\dfrac{1}{2}\dfrac{\phi_{p+1,q+1}-\phi_{p+1,q}+\phi_{p,q+1}-\phi_{p,q}}{\Delta y}
	\end{pmatrix}.
	\label{eq:compatible_grad}
\end{equation}
Once the compatible gradient operator has been presented, we will define the \textit{discrete curl operator} $\nabla^h\times$ of a vector field $\nabla^h \phi$ at the center of the cells, making use of the discrete gradients in the surrounding cells. The component in the $z-$direction is given by 
\begin{align}
	\left(\nabla^h\times \nabla^h \phi\right)_{p,q} \!\!\!\!\cdot\boldsymbol{e}_3 
	&=
	\dfrac{\left(\partial_y^h \phi\right)_{p+\frac{1}{2},q+\frac{1}{2}}\!\!-\left(\partial_y^h \phi\right)_{p-\frac{1}{2},q+\frac{1}{2}}}{2\Delta x}+\dfrac{\left(\partial_y^h \phi\right)_{p+\frac{1}{2},q-\frac{1}{2}}\!\!-\left(\partial_y^h \phi\right)_{p-\frac{1}{2},q-\frac{1}{2}}}{2\Delta x}\nonumber\\
	&-\dfrac{\left(\partial_x^h \phi\right)_{p+\frac{1}{2},q+\frac{1}{2}}\!\!-\left(\partial_x^h \phi\right)_{p+\frac{1}{2},q-\frac{1}{2}}}{2\Delta y}-\dfrac{\left(\partial_x^h \phi\right)_{p-\frac{1}{2},q+\frac{1}{2}}\!\!-\left(\partial_x^h \phi\right)_{p-\frac{1}{2},q-\frac{1}{2}}}{2\Delta y},
	\label{eq:compatible_curl}
\end{align}
with $\boldsymbol{e}_3=(0,0,1)$. 
Once the compatible discrete gradient and curl operators have been defined, it is necessary to verify that the continuous identity
\begin{equation}
	\nabla \times \nabla \phi = 0,
	\label{eq:continuous_curlgrad}
\end{equation}
is also satisfied at the discrete level. Combining~\eqref{eq:compatible_grad} and~\eqref{eq:compatible_curl}, it is easy to prove that for an arbitrary scalar field $\phi$ defined in the barycenter of the main grid one has
\begin{align*}
	\left(\nabla^h\times\right.&\left.\nabla^h \phi\right)_{p,q}\cdot \boldsymbol{e}_3 \nonumber\\	
	=\dfrac{1}{4}&\left( \dfrac{ \left(\phi_{p+1,q+1}-\phi_{p+1,q}+\phi_{p,q+1}-\phi_{p,q}\right)
		-\left(\phi_{p,q+1}-\phi_{p,q}+\phi_{p-1,q+1}-\phi_{p-1,q}\right) }{\Delta x \Delta y}\right.\nonumber\\
	&\left.+\dfrac{ \left(\phi_{p+1,q}-\phi_{p+1,q-1}+\phi_{p,q}-\phi_{p,q-1}\right)
		-\left(\phi_{p,q}-\phi_{p,q-1}+\phi_{p-1,q}-\phi_{p-1,q-1}\right) }{\Delta x \Delta y}\right. \nonumber\\                       	
	&\left.-\dfrac{ \left(\phi_{p+1,q+1}-\phi_{p,q+1}+\phi_{p+1,q}-\phi_{p,q}\right)
		-\left(\phi_{p+1,q}-\phi_{p,q}+\phi_{p+1,q-1}-\phi_{p,q-1}\right) }{\Delta x \Delta y}\right.\nonumber\\
	&\left.-\dfrac{ \left(\phi_{p,q+1}-\phi_{p-1,q+1}+\phi_{p,q}-\phi_{p-1,q}\right) 
		-\left(\phi_{p,q}-\phi_{p-1,q}+\phi_{p,q-1}-\phi_{p-1,q-1}\right) }{\Delta x\Delta y}\right)=0,
	\label{eq:proof_discrete_curlgrad}
\end{align*}
that is, \eqref{eq:continuous_curlgrad} is satisfied also at the discrete level or equivalently,
\begin{equation}
	\nabla^h \times \nabla^h \phi = 0,
	\label{eq:discrete_curlgrad}
\end{equation}
for all cells of the computational domain, which is the discrete curl-grad compatibility. This equality shows that any gradient field defined using the discrete compatible operator \eqref{eq:compatible_grad} is exactly curl-free for the discrete compatible operator \eqref{eq:compatible_curl}. 


\subsection{Discretization of the relative velocity with the compatible operators}
The relative velocity equation will be discretized making use of the compatible operators defined before. Thus, Equation~\eqref{eq:relative_vel_ini} reads as	
\begin{align}
	\left(\w{k}\right)_{p+\frac{1}{2},q+\frac{1}{2}}^{n+1}&=\left(\w{k}\right)_{p+\frac{1}{2},q+\frac{1}{2}}^{n}-\Delta t \,\partial_k^h \left( \w{l} \u{l}+\varphi\right)_{p+\frac{1}{2},q+\frac{1}{2}}^n\nonumber\\
	&- \frac{\Delta t}{4}\sum_{r=0}^1 \sum_{s=0}^1 \left(\u{l}\right)_{p+r,q+s}^{n}\left(\partial_l^{h} \left(\w{k}\right)_{p+\frac{1}{2},q+\frac{1}{2}}^{n}-\partial_k^{h}\left(\w{l}\right)_{p+\frac{1}{2},q+\frac{1}{2}}^{n}\right),
	\label{eq:compatible_relative_vel}
\end{align}
for each component $k$ of the gradient field $\bw$.	In the right-hand side of the equation~\eqref{eq:compatible_relative_vel}, the last term can be computed using the compatible gradient~\eqref{eq:compatible_grad}. Now, we will prove that for a curl-free vector $\bw$ that satisfies $\nabla^h\times \bw^{h,n}=0$, it is also satisfied $\nabla^h\times \bw^{h,n+1}=0$.
If the discrete curl operator $\nabla^h\times$ is applied to~\eqref{eq:compatible_relative_vel}, taking into account that $\nabla^h\times\bw^{h,n}=0$, it results
\begin{equation}
	\nabla^h\times\left(\w{k}\right)_{p+\frac{1}{2},q+\frac{1}{2}}^{n+1}=
	-\Delta t \, \nabla^h\times\left(\partial_k^h \left( \w{l} \u{l}+\varphi\right)_{p+\frac{1}{2},q+\frac{1}{2}}^n\right).
	\label{eq:compatible_relative_vel2}
\end{equation}
Since the right hand side of Equation~\eqref{eq:compatible_relative_vel2} is zero due to the definition of the compatible discrete curl-grad~\eqref{eq:discrete_curlgrad}, then $\nabla^h\times\left(\w{k}\right)^{h,n+1}=0$. That is, if the vector field $\bw$ is curl-free at the initial time, it will be curl-free for every instant of time.

The scheme just presented is based on a central discretization. We will define a compatible artificial viscosity to suppress instabilities that may appear and to guarantee the stability of the method, while keeping the discretely curl-free structure. Recall that the vector Laplacian of $\bw$ at the continuous level can be written as
\begin{equation}
	\nabla^2 \bw = \nabla \left(\nabla \cdot \bw\right) - \nabla \times \nabla \times \bw
	\label{eq:vector_Laplacian_cont}
\end{equation}
To discretize~\eqref{eq:vector_Laplacian_cont}, defining a discrete divergence operator is necessary. We will define the \textit{discrete divergence operator} $\nabla^h\cdot$ acting on a discrete vector field $\bw^h=(\w{x}^h,\w{y}^h)^T$ at the corners, considering as a stencil a piecewise linear reconstruction of $\bw$ at the cell centers. We will calculate, therefore, the extrapolated values of the cell centers of $\bw$, for each $(\x{p},\y{q})$. Using such values, the discrete divergence operator at vertex $(\x{p+\frac{1}{2}},\y{q+\frac{1}{2}})$ is then defined as:
\begin{align}
	\nabla^h\cdot \bw^h =\left(\partial^h_k \w{k}\right)_{p+\frac{1}{2},q+\frac{1}{2}}
	&\!=\dfrac{1}{2}\left(\dfrac{(\w{x})_{p+1,q+1} - (\w{x})_{p,q+1} }{\Delta x} + \dfrac{(\w{x})_{p+1,q} - (\w{x})_{p,q} }{\Delta x}\right)\nonumber\\
	&\!+\dfrac{1}{2}\left(\dfrac{(\w{y})_{p+1,q+1} - (\w{y})_{p+1,q} }{\Delta y} + \dfrac{(\w{y})_{p,q+1} - (\w{y})_{p,q} }{\Delta y}\right).
	\label{eq:compatible_div}
\end{align}
Having defined the discrete divergence operator~\eqref{eq:compatible_div}, the discrete version of~\eqref{eq:vector_Laplacian_cont} at $(\x{p+\frac{1}{2}},\y{q+\frac{1}{2}})$ is given as follows
\begin{equation}
	\left(\nabla^2_h \bw\right)_{p+\frac{1}{2},q+\frac{1}{2}} = \left(\nabla^h \left(\nabla^h \cdot \bw\right)\right)_{p+\frac{1}{2},q+\frac{1}{2}} - \left(\nabla^h \times \nabla^h \times \bw\right)_{p+\frac{1}{2},q+\frac{1}{2}}
	\label{eq:vector_Laplacian_discrete}
\end{equation}
Multiplying \eqref{eq:vector_Laplacian_discrete} by the mesh size $h=\max(\Delta x,\Delta y)$ and by an appropriate constant signal speed $c_h$, the expression~\eqref{eq:compatible_relative_vel} turns out to be
\begin{align*}
	\left(\w{k}\right)_{p+\frac{1}{2},q+\frac{1}{2}}^{n+1}&=\left(\w{k}\right)_{p+\frac{1}{2},q+\frac{1}{2}}^{n}-\Delta t \,\left(\partial_k^h \left( \w{l} \u{l}+\varphi\right)-hc_h\partial^h_l \w{l}\right)_{p+\frac{1}{2},q+\frac{1}{2}}^n\nonumber\\
	&- \frac{\Delta t}{4}\sum_{r=0}^1 \sum_{s=0}^1 \left(\u{l}\right)_{p+r,q+s}^{n}\left(\partial_l^{h} \left(\w{k}\right)_{p+\frac{1}{2},q+\frac{1}{2}}^{n}-\partial_k^{h}\left(\w{l}\right)_{p+\frac{1}{2},q+\frac{1}{2}}^{n}\right).
\end{align*}
This compatible discretization can be used to discretize equation \eqref{eq:relative_vel_ini}, and thus, as long as it is curl-free at the initial time, $\bw$ remains curl-free for all times.

\subsection{Discretization of the remaining terms}
Once equation \eqref{eq:flux_splittingv} has been discretized with the compatible operators described in the previous sections, to discretize the remaining terms (Equations~\eqref{eq:flux_splittingb})  we will make use of a classical second-order MUSCL-Hancock scheme. Below, we will briefly describe the numerical scheme. It reads 
\begin{align*}
	\bQ_{p,q}^{n+1} = \bQ_{p,q}^{n} 
	&	- \frac{\Delta t}{\Delta x} \left(  \left(\bF_x\right)_{p+\frac{1}{2},q}-\left(\bF_x\right)_{p-\frac{1}{2},q} \right)  
	- \frac{\Delta t}{\Delta y} \left(  \left(\bF_y\right)_{p,q+\frac{1}{2}}-\left(\bF_y\right)_{p,q-\frac{1}{2}} \right) \\		
	&	- \frac{\Delta t}{\Delta x} \left(  \left(\bDb_x\right)_{p+\frac{1}{2},q}+\left(\bDb_x\right)_{p-\frac{1}{2},q} \right)  
	- \frac{\Delta t}{\Delta y} \left(  \left(\bDb_y\right)_{p,q+\frac{1}{2}}+\left(\bDb_y\right)_{p,q-\frac{1}{2}} \right) \\		
	&   - \Delta t \,\bBb(\bQ_{p,q}^{n+\frac{1}{2}})\cdot\nabla\bQ_{p,q}^{n+\frac{1}{2}} +\Delta t\,\bSb(\bQ_{p,q}^{n+\frac{1}{2}}),		
\end{align*}
where $\left(\bF_x\right)_{p+\frac{1}{2},q}$ and $\left(\bF_y\right)_{p,q+\frac{1}{2}}$ are the numerical fluxes in $x$ and $y$ direction, respectively, including $\bFb$ and $\bFv$. The path-conservative jump terms according to Castro and Par\'es  \cite{Castro2006,Pares2006,Castro2008} are given in the $x$ and $y$ direction by 
\begin{equation*}
	\bDb_{x,y}(\bQ_h^{-},\bQ_h^{+}) = \frac{1}{2}\widetilde{\bB{}}_{x,y}    (\bQ_h^{+}-\bQ_h^{-}), \quad \text{with} \quad \widetilde{\bB{}}_{x,y}  = \int_0^1 \bBb{}(	\Psi(s, \bQ_h^+, \bQ_h^-)) \cdot \bn_{x,y} \ \ddx{S}.
\end{equation*}
They are defined as a function of two generic left and right boundary-extrapolated values $\bQ_h^-$ and $\bQ_h^+$, respectively. These terms take into account jumps of $\bQ$ at the element boundaries. In this paper the simple straight-line segment path is chosen
\begin{equation*}
	\Psi(s, \bQ_h^+, \bQ_h^-) =  \bQ_h^-+s( \bQ_h^+- \bQ_h^+), \quad s\in[0,1].
\end{equation*}
In order to compute the numerical fluxes $\bF_x$ and $\bF_y$, we will use the Rusanov flux, see~\cite{Rusanov1962Nonstationary}:	
\begin{align*}
	\left(\bF_x\right)_{p+\frac{1}{2},q}
	&= \frac{1}{2} \left( 
	\bFb_x \left(\bQ_{p+\frac{1}{2},q}^{n+\frac{1}{2},\,-}\right)
	+ \bFb_x \left(\bQ_{p+\frac{1}{2},q}^{n+\frac{1}{2},\,+}\right) 
	+ \bFv_x \left(\bQ_{p+\frac{1}{2},q+\frac{1}{2}}^n\right) 
	+ \bFv_x \left(\bQ_{p+\frac{1}{2},q-\frac{1}{2}}^n\right) \right) \nonumber\\
	&- \frac{1}{2} s_{\max}^x \left(\bQ_{p+\frac{1}{2},q}^{n+\frac{1}{2},\,+} - \bQ_{p+\frac{1}{2},q}^{n+\frac{1}{2},\,-}\right),\\
	\left(\bF_y\right)_{p,q+\frac{1}{2}}
	&= \frac{1}{2} \left( 
	\bFb_y \left(\bQ_{p,q+\frac{1}{2}}^{n+\frac{1}{2},\,-}\right) 
	+ \bFb_y \left(\bQ_{p,q+\frac{1}{2}}^{n+\frac{1}{2},\,+}\right) 
	+ \bFv_y \left(\bQ_{p+\frac{1}{2},q+\frac{1}{2}}^n\right) 
	+ \bFv_y \left(\bQ_{p-\frac{1}{2},q+\frac{1}{2}}^n\right) \right) \nonumber\\
	&- \frac{1}{2} s_{\max}^y \left(\bQ_{p,q+\frac{1}{2}}^{n+\frac{1}{2},\,+} - \bQ_{p,q+\frac{1}{2}}^{n+\frac{1}{2},\,-}\right),
	\label{eq:rusanov_flux} 
\end{align*} 
where $s_{\max}^x$ and $s_{\max}^y$ are the maximum wave speeds in the $x-$ and $y$ direction, respectively, computed as the maximum absolute value of the eigenvalues in each direction. The values $\bQ^{\pm}$ are the boundary extrapolated values, given by
\begin{align*}
	\bQ_{p+\frac{1}{2},q}^{n+\frac{1}{2},\,-} &= \bQ_{p,q}^n +\frac{1}{2}\Delta x \partial_x \bQ_{p,q}^n + \frac{1}{2}\Delta t \partial_t \bQ_{p,q}^n,\\
	\bQ_{p+\frac{1}{2},q}^{n+\frac{1}{2},\,+} &= \bQ_{p+1,q}^n -\frac{1}{2}\Delta x \partial_x \bQ_{p+1,q}^n + \frac{1}{2}\Delta t \partial_t \bQ_{p+1,q}^n,\\
	\bQ_{p,q+\frac{1}{2}}^{n+\frac{1}{2},\,-} &= \bQ_{p,q}^n +\frac{1}{2}\Delta y \partial_y \bQ_{p,q}^n + \frac{1}{2}\Delta t \partial_t \bQ_{p,q}^n,\\
	\bQ_{p,q+\frac{1}{2}}^{n+\frac{1}{2},\,+} &= \bQ_{p,q+1}^n -\frac{1}{2}\Delta y \partial_y \bQ_{p,q+1}^n + \frac{1}{2}\Delta t \partial_t \bQ_{p,q+1}^n,\\
\end{align*}
with $\partial_x \bQ_{p,q}^n$ and $\partial_y \bQ_{p,q}^n$ the slopes defined as
\begin{align*}
	\partial_x \bQ_{p,q}^n&=\mathrm{minmod}\left(\dfrac{\bQ_{p+1,q}^n-\bQ_{p,q}^n}{\Delta x}, \dfrac{\bQ_{p,q}^n-\bQ_{p-1,q}^n}{\Delta x}\right),\\
	\partial_y \bQ_{p,q}^n&=\mathrm{minmod}\left(\dfrac{\bQ_{p,q+1}^n-\bQ_{p,q}^n}{\Delta y}, \dfrac{\bQ_{p,q}^n-\bQ_{p,q-1}^n}{\Delta y}\right),
\end{align*}
and with the time derivative $\partial_t \bQ_{p,q}^n$ computed as follows
\begin{align*}
	\partial_t \bQ_{p,q}^n = 
	&-\dfrac{\bFb_x \left( \bQ_{p,q}^n +\frac{1}{2}\Delta x \partial_x \bQ_{p,q}^n\right) - \bFb_x \left( \bQ_{p,q}^n -\frac{1}{2}\Delta x \partial_x \bQ_{p,q}^n\right)}{\Delta x}\nonumber\\
	&-\dfrac{\bFb_y \left( \bQ_{p,q}^n +\frac{1}{2}\Delta y \partial_y \bQ_{p,q}^n\right) - \bFb_y \left( \bQ_{p,q}^n -\frac{1}{2}\Delta y \partial_y \bQ_{p,q}^n\right)}{\Delta y}\nonumber\\
	&-\dfrac{\bFv_x \left( \bQ_{p+\frac{1}{2},q+\frac{1}{2}}^n\right)
		\!+\! \bFv_x \left( \bQ_{p+\frac{1}{2},q-\frac{1}{2}}^n\right)
		\!-\! \bFv_x \left( \bQ_{p-\frac{1}{2},q+\frac{1}{2}}^n\right)
		\!-\! \bFv_x \left( \bQ_{p-\frac{1}{2},q-\frac{1}{2}}^n\right)}{2\Delta x}\nonumber\\
	&-\dfrac{\bFv_y \left( \bQ_{p+\frac{1}{2},q+\frac{1}{2}}^n\right)
		\!+\! \bFv_y \left( \bQ_{p-\frac{1}{2},q+\frac{1}{2}}^n\right)
		\!-\! \bFv_y \left( \bQ_{p+\frac{1}{2},q-\frac{1}{2}}^n\right)
		\!-\! \bFv_y \left( \bQ_{p-\frac{1}{2},q-\frac{1}{2}}^n\right)}{2\Delta y}\nonumber\\
	&+\bSb(\bQ_{p,q}^n)-\bBb(\bQ_{p,q}^n)\cdot\nabla\bQ_{p,q}^{n}.
\end{align*}
This completes the description of the numerical method. To summarize: in our new structure-preserving algorithm we employ a path-conservative MUSCL-Hancock-type finite volume scheme for the evolution of all quantities, apart from the relative velocity field, which is updated using a special discretization on suitably staggered meshes which allows to preserve the curl-free property of the relative velocity field exactly also at the discrete level. 

\clearpage 

\section{Numerical results}
\label{sec:Numerical_Results}
In this section, we will show some test cases that demonstrate the performance of the proposed method, the convergence order and illustrate that the method is exactly curl-free. First, we will use the proposed methodology to solve a 1D Riemann problem, perform some simulations to show the experimental order of convergence (EOC), and solve a 2D circular explosion problem. In all cases, reference solutions are used for a fair comparison with the approximate solution. Then, we will simulate a dam break and compare the results with those obtained using the reduced Baer--Nunziato model. Finally, we will perform a simulation to show the Kelvin--Helmholtz instabilities and qualitatively compare the results with those existing in the literature to show the good behavior of the proposed method. In all the results, the gravity $\bg$ is set to $\bm{0}$, except in the dam break test, where it has to be considered.

\subsection{1D Riemann problem}
\begin{figure}[ht!]
	\begin{tabular}{lll} 
		\includegraphics[trim=15 15 60 50,clip,width=0.42\textwidth]{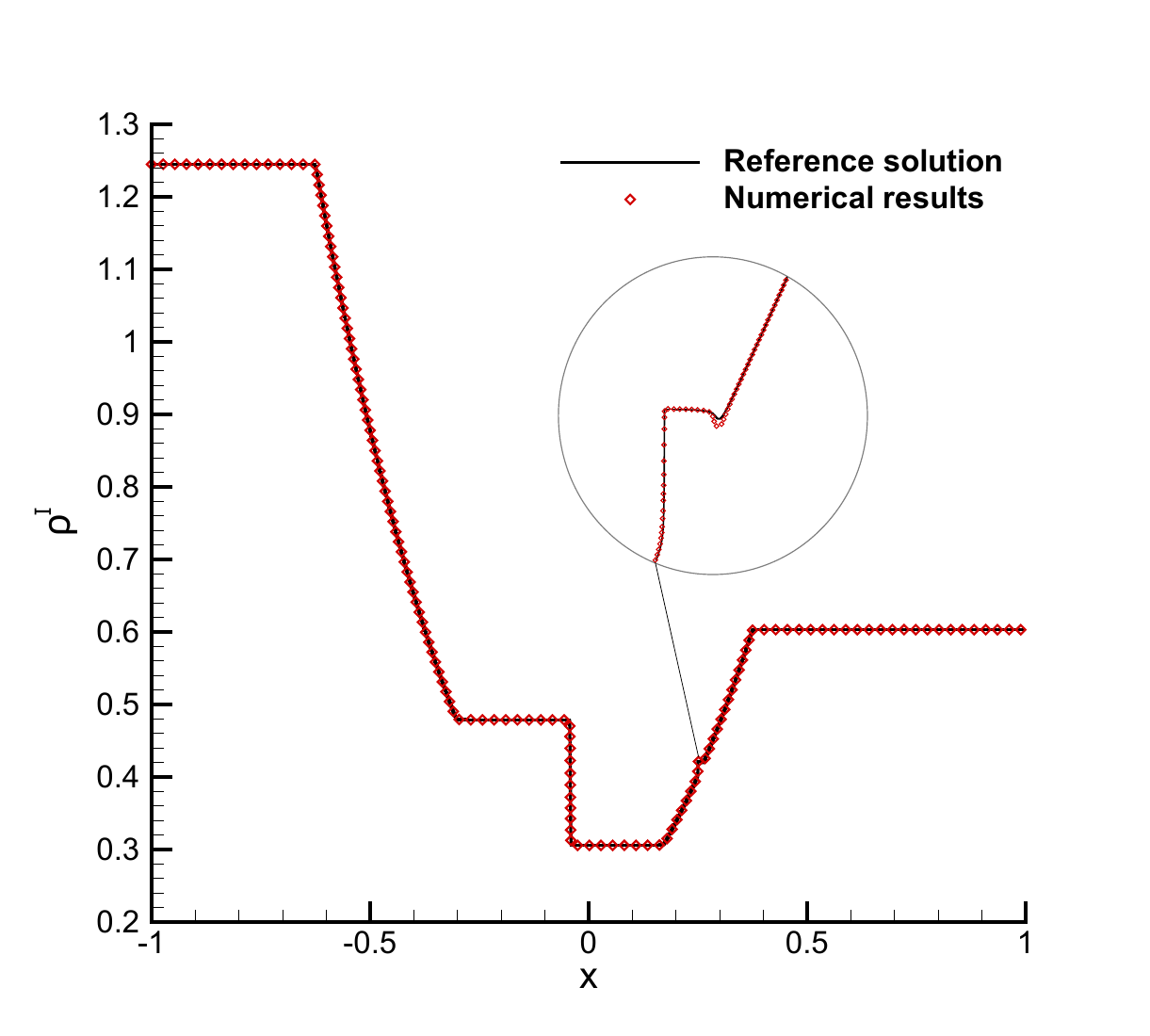}   & 
		\includegraphics[trim=15 15 60 50,clip,width=0.42\textwidth]{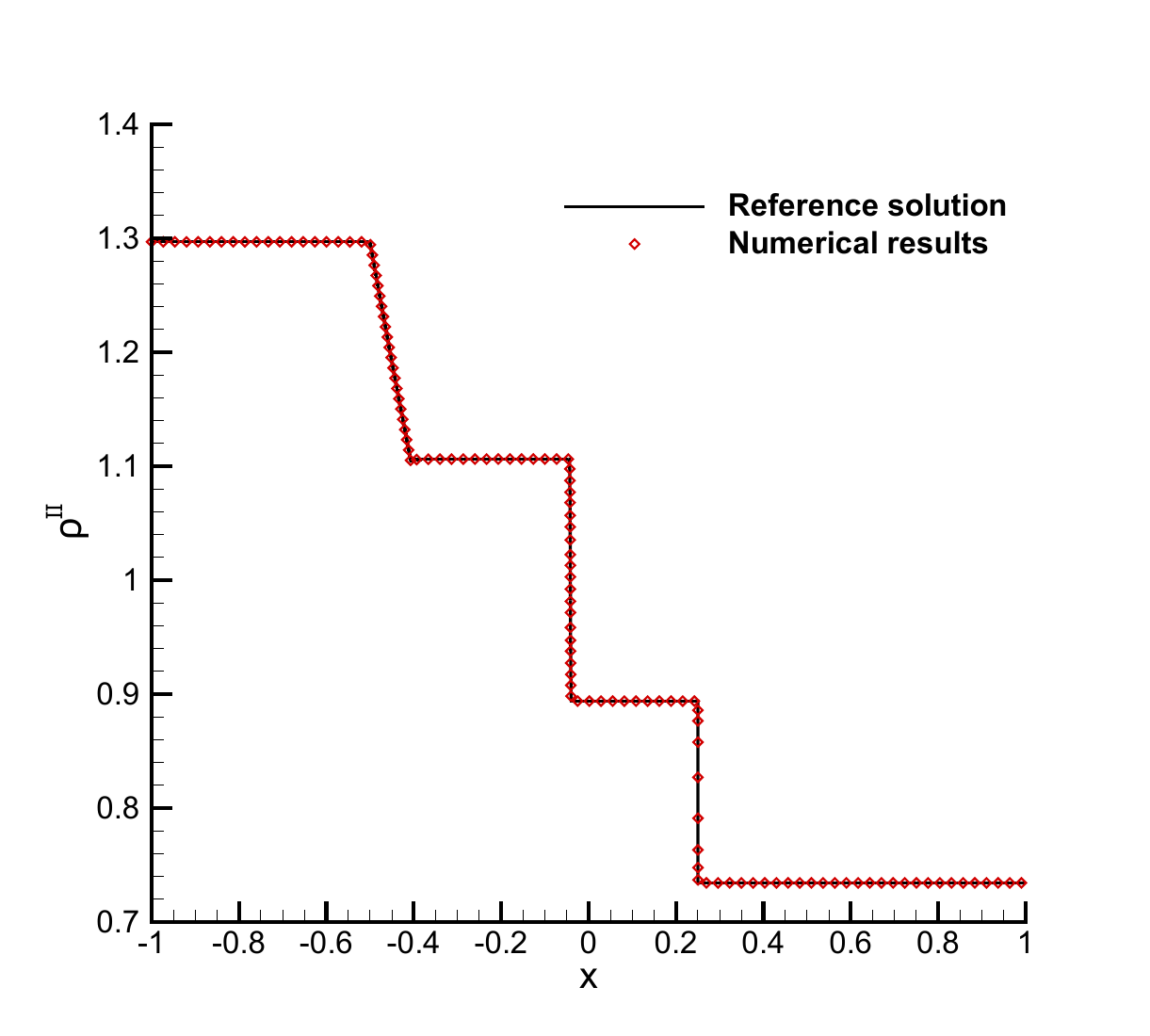}   \\ 
		\includegraphics[trim=15 15 60 50,clip,width=0.42\textwidth]{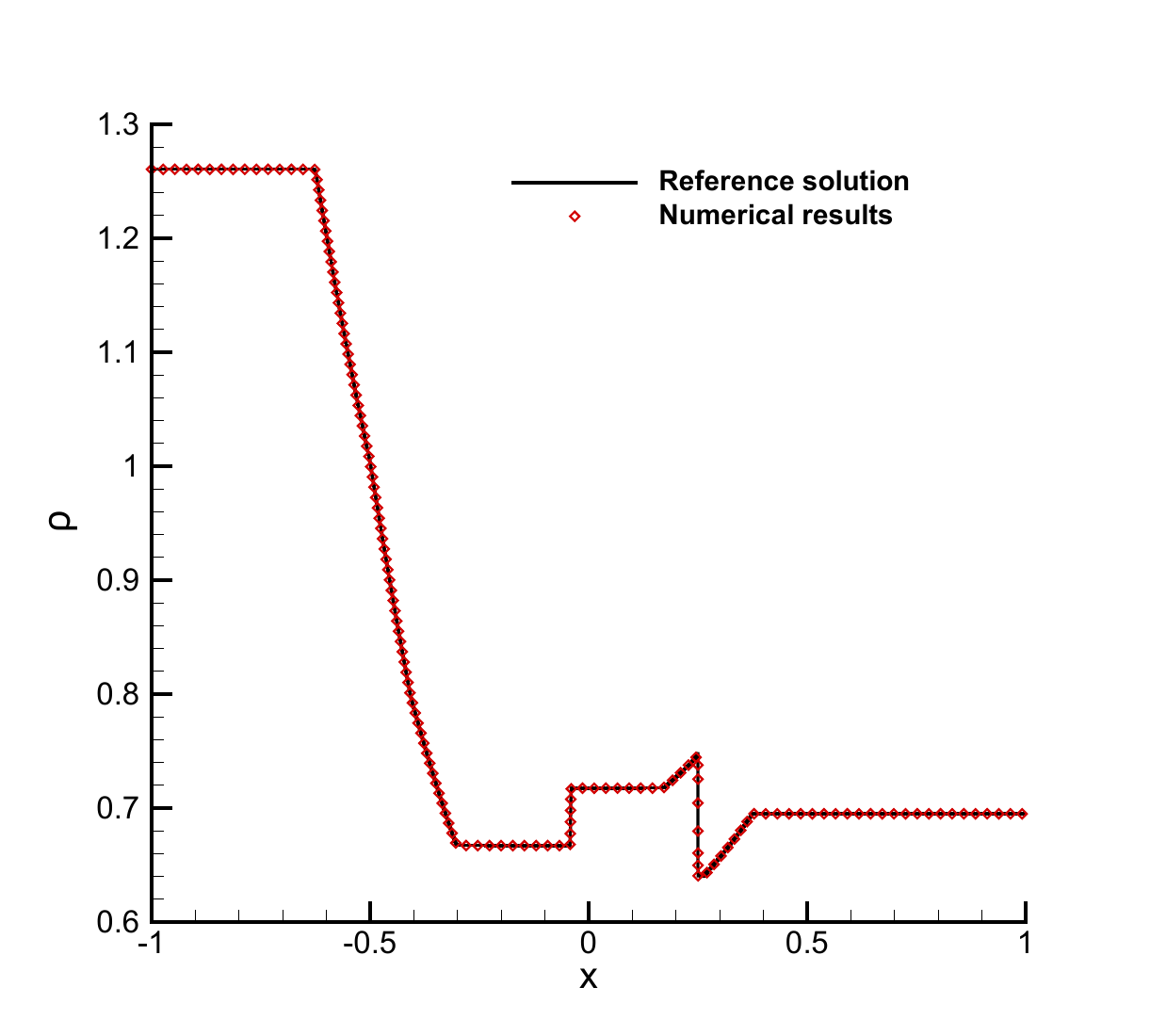}    &   
		\includegraphics[trim=15 15 60 50,clip,width=0.42\textwidth]{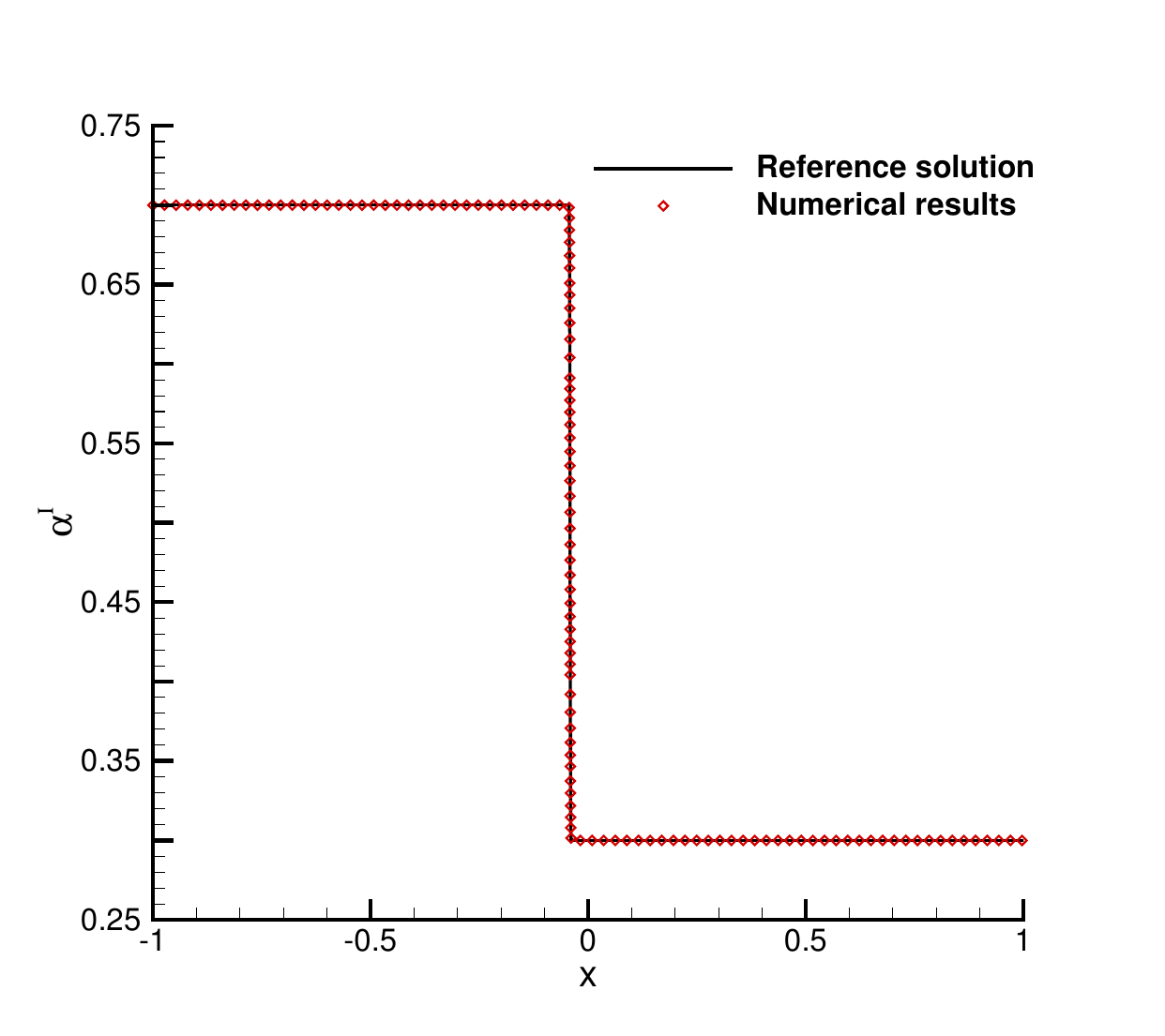}  \\
		\includegraphics[trim=15 15 60 50,clip,width=0.42\textwidth]{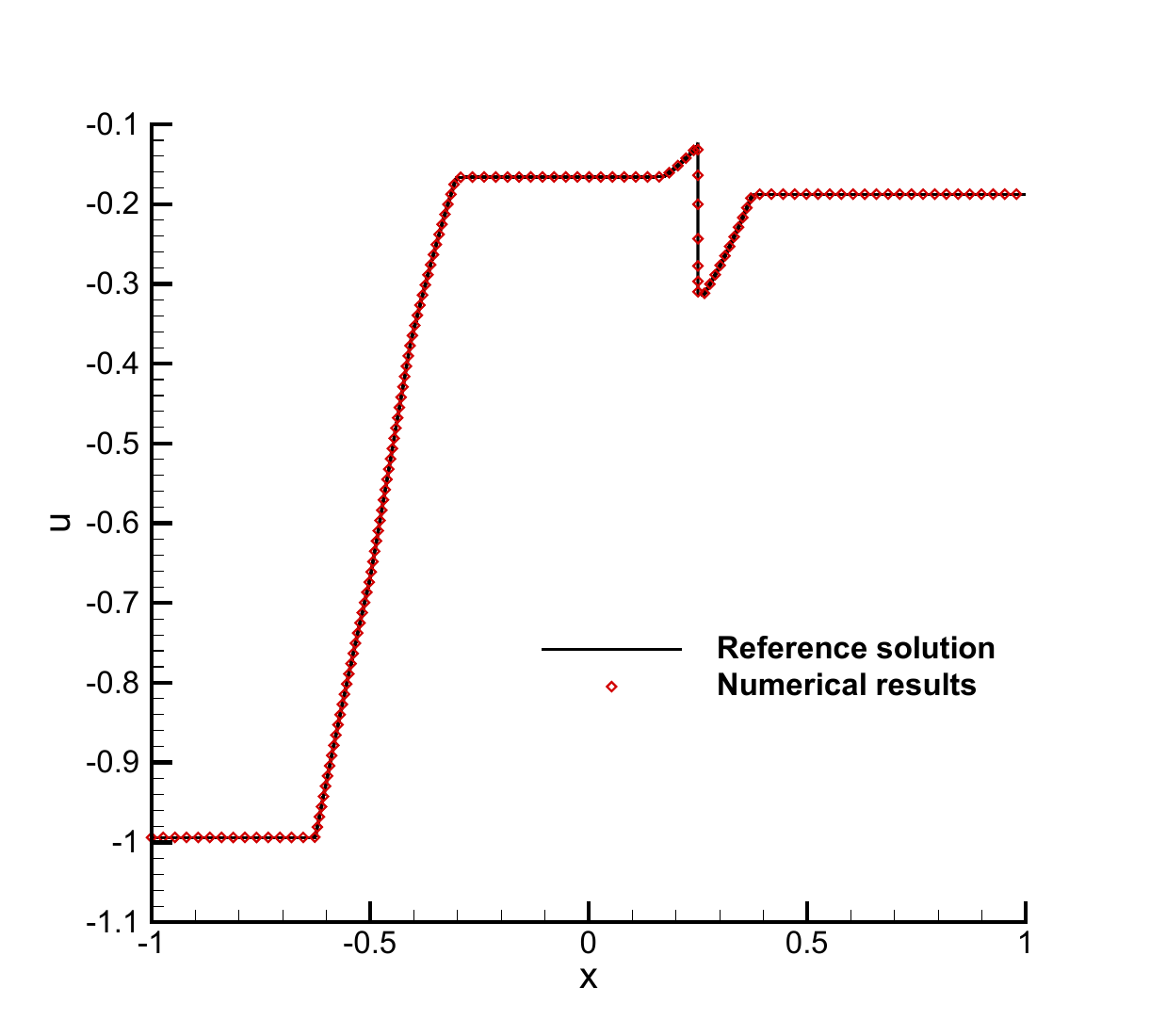}      &
		\includegraphics[trim=15 15 60 50,clip,width=0.42\textwidth]{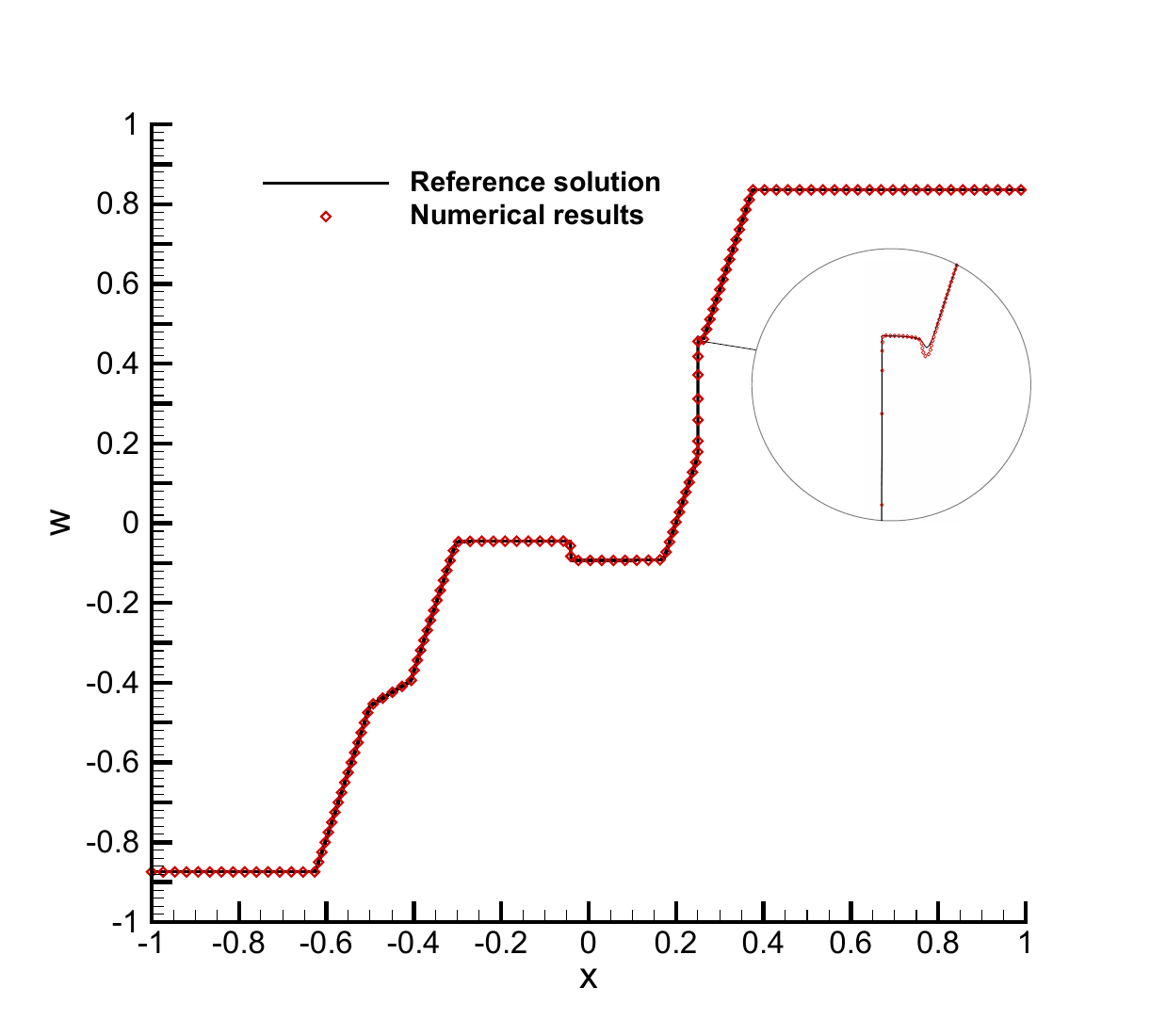}      
	\end{tabular}
	\caption{1D Riemann problem solved with the structured preserving finite volume scheme on a Cartesian staggered mesh with 30000 cells at final time $t=0.25$. Top: densities of each phase, $\rhoI$ and $\rhoII$. Center: mixture density $\rho$ and volume fraction $\alphaI$. Bottom: mixture velocity $\u{}$ and relative velocity $\w{}=\uI{}-\uII{}$. The zooms show the shock inside the rarefaction which are a very special feature of the exact solution of this Riemann problem, see \cite{Thein2022}.} 
	\label{fig:RP1D}
\end{figure}

The first test case under consideration is a Riemann problem in 1D, where a shock in one phase appears within a rarefaction wave of the other phase, see \cite{Thein2022} for the exact solution of the Riemann problem and a detailed discussion of possible wave patterns. The computational domain is $\Omega = [-1,1]$ and has been discretized using a mesh with $30000$ cells in order to show sufficiently converged numerical results. The final time is $t=0.25$, and the CFL has been set to $0.25$. In both phases, we will consider the EOS for an ideal gas~\eqref{eq:energy_ideal_gas} with $\gammaI=1.4$, and $\gammaII=2$. The left and right states used are shown in Table~\ref{tab:RP1} (see~\cite{Rio-Martin2023} and~\cite{Thein2022} for more details on the test). A comparison has been made between the numerical results obtained with the method proposed in this work and the reference solution computed with a second-order MUSCL--Hancock scheme based on the Rusanov flux as an approximate Riemann solver and using a mesh spacing of $\Delta x = 2\cdot 10^{-6}$ (see~\cite{Thein2022} for details about the reference solution). The results, for the densities of each phase $\rhoI,\,\rhoII$, the mixture density $\rho$, the volume fraction $\alphaI$, the mixture velocity $\u{}$ and the relative velocity $\w{}$ are shown in Figure~\ref{fig:RP1D}. Excellent agreement between the calculated and the reference solutions is observed in all cases. The choice of the number of points in this one-dimensional test is due to the need to capture the resonance phenomenon of the shock inside a rarefaction wave accurately enough, as illustrated in the zooms of the first phase density and relative velocity. 
\begin{table}[ht!]
	\centering
	\caption{Left and right states of the one dimensional Riemann problem}
	\begin{tabular}{cccccccccc}
		\toprule
		& $\alphaI$ & $\rhoI$ & $\rhoII$ & $\uI{}$ & $\uII{}$  \\	\midrule 
		$Q_L$ & 0.7 & 1.2449 & 1.2969 & -1.2638 & -0.38947 \\
		$Q_R$ & 0.3 & 0.60312 & 0.73436 & 0.43059 & -0.40507 \\\bottomrule 
	\end{tabular}
	\label{tab:RP1}
\end{table} 

\subsection{Stationary vortex solution}
\begin{figure}[h!]
	\begin{center} 
		\includegraphics[trim=5 15 50 50,clip,width=0.48\textwidth]{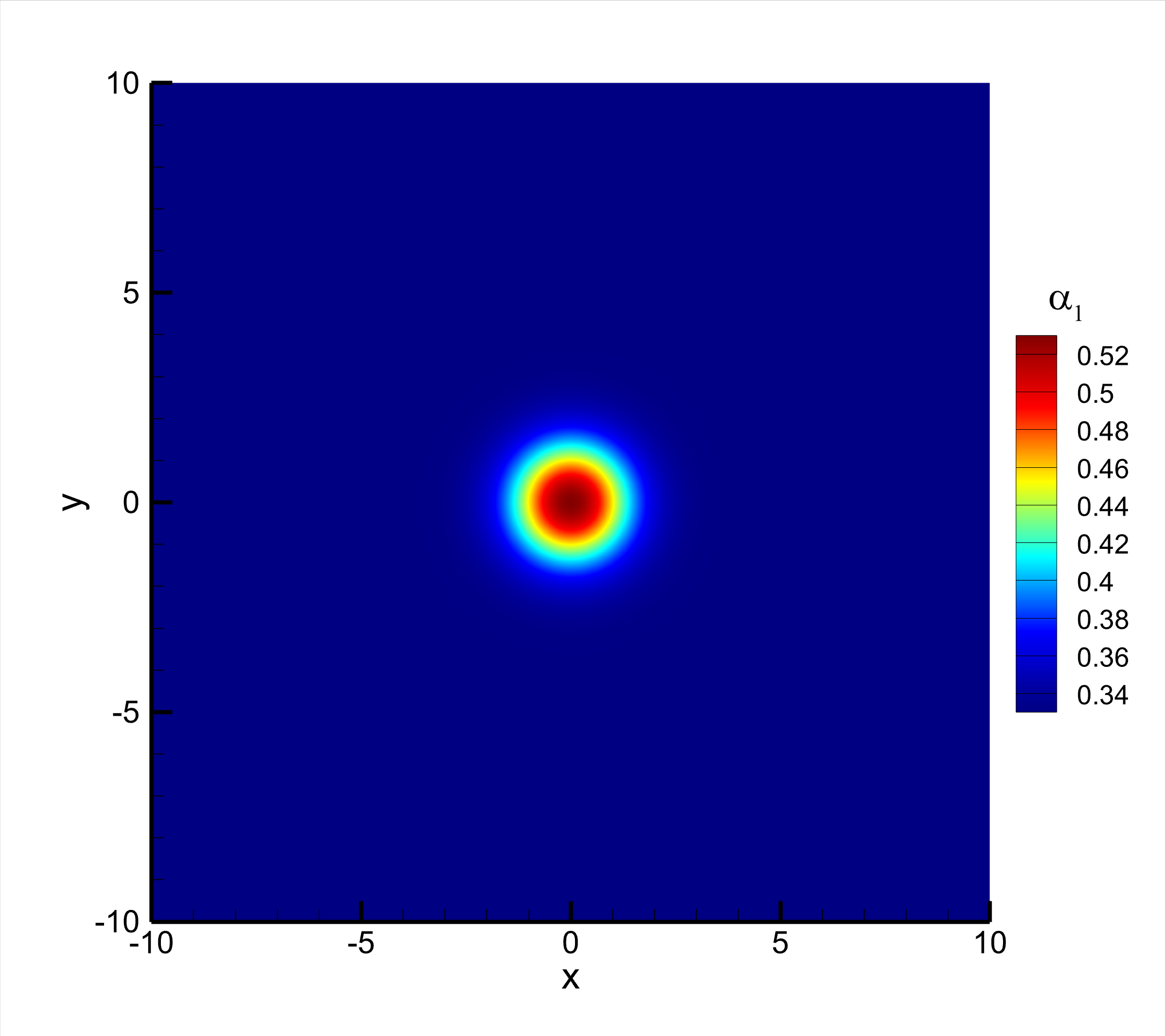}    
		\includegraphics[trim=5 15 50 50,clip,width=0.48\textwidth]{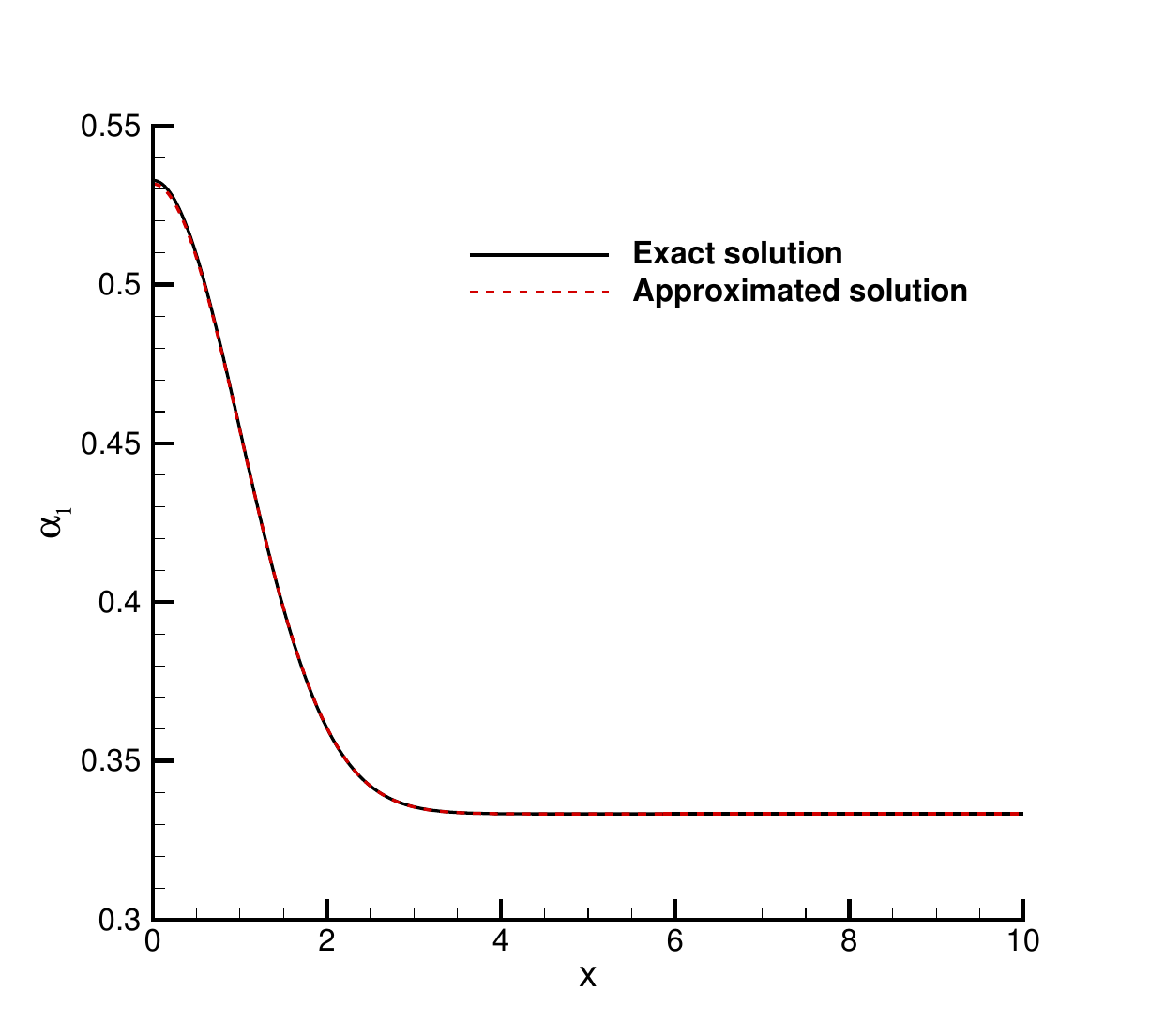}   
		\caption{Left: Solution of the stationary vortex with the curl-free scheme solved on a staggered mesh with $1024\times1024$ cells at time $t=1000$. Right: Comparison between the approximated solution (dashed red line) and the exact solution (solid black line) at time $t=1000$.} 
		\label{fig:vortex_t1000}
	\end{center}
\end{figure}

\begin{figure}[h!]
	\begin{center}
		\begin{tabular}{cc} 
			\includegraphics[trim=1 1 1 1,clip,width=0.65\textwidth]{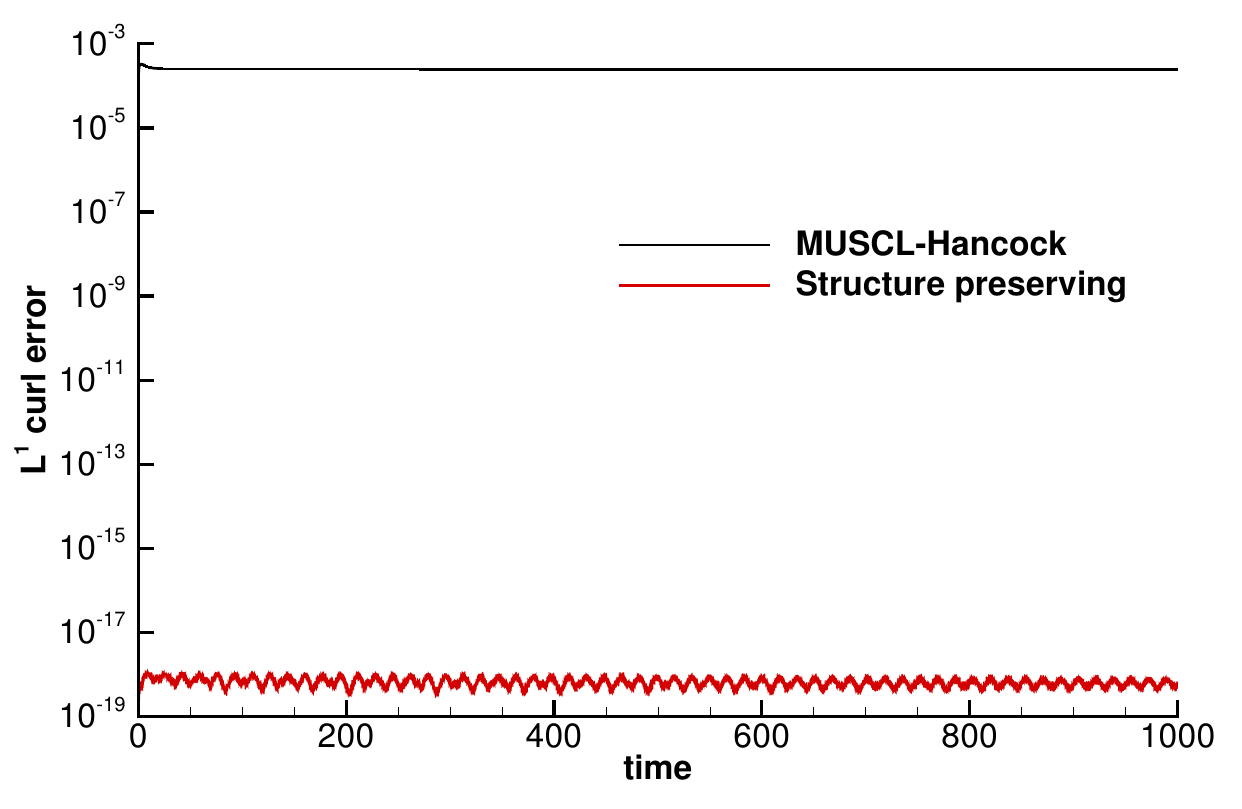}   
		\end{tabular} 
		\caption{Time-evolution of the $L^1$ norm of the discrete curl errors using the staggered compatible curl-free discretization (red line) and without the compatible curl-free discretization (black line) for the 2D stationary vortex problem, using a mesh with $1024\times 1024$ cells.} 
		\label{fig:curlErrors_vortex}
	\end{center}
\end{figure}

In order to have a quantitative assessment of the scheme's performance, we consider here a two-dimensional stationary solution to check the experimental order of convergence and evaluate the curl-free property. For this, we consider the rotationally symmetric exact solution of system~\eqref{eq:TwoFluid} presented in \cite{Rio-Martin2023}. In polar coordinates $(r,\theta)$, the solution can be expressed in terms of primitive variables as follows 
\begin{align*}
	\alphaI(r) &= \frac{1}{3} + \frac{e^{-\frac{r^2}{2}}}{2\sqrt{2\pi}},\qquad \rhoI(r) = \rhoII(r) = \prn{1-\frac{e^{1-r^2}}{4}}^{5/7}\!\!\!\!\!, \\
	\uI{\theta}(r) &= \uII{\theta}(r) = 2^{3/14}\!\!\sqrt{\frac{r^2 e^{1-r^2}}{\prn{4-e^{1-r^2}}^{5/7}}}.
\end{align*}
For this test case, we use the same ideal gas equation of state for both fluids with $\gammaI=\gammaII=1.4$. The acceleration of gravity is neglected here. The computational domain is $\Omega=[-10,10]\times[-10,10]$. Periodic boundary conditions are used in both directions. In Figure~\ref{fig:vortex_t1000}, the left plot shows the solution of the stationary vortex computed using the structured preserving FV scheme for a mesh resolution of $1024\times1024$ cells at time $t=1000$. On the right plot, we show a graphical representation of the obtained solution using the structure-preserving FV scheme, plotted with dashed red lines, compared with the exact solution, plotted with solid black lines.
Figure~\ref{fig:curlErrors_vortex} shows the $L^1$ norm of the curl errors for the simulation with final time $t=1000$ using a mesh resolution of $1024\times1024$ cells for the solution using the new structure-preserving scheme described in this work and without the compatible curl-free discretization. The obtained results clearly show that the new method is able to maintain the curl-free property of the relative velocity field up to machine precision. The number of time steps performed is around 400000.
A convergence table for this test case on different mesh resolutions is given in Table~\ref{tab:conv}, which shows that the second order is well recovered in all the representative variables.
\begin{table}[h!]
	\begin{center}
	\caption{Convergence table for the vortex test case with the proposed curl-free scheme. $\epsilon$ represents the $L^2$ error norm of the corresponding variable and $O$ is the experimental convergence order. All model parameters are taken similarly as for Figure~\ref{fig:vortex_t1000}, except for the mesh sizes, and the final time set to $t=1$.}
	\small
	\begin{tabular}{ccccccc}
		\toprule
		$N_x=N_y$ & $\epsilon_{\alpha_1}$ & $\epsilon_{\rho_1}$ & $\epsilon_{\rho_2}$ & $\epsilon_{u_1}$ & $\epsilon_{u_2}$ \\
		\midrule
		192 & $1.92 \times 10^{-4}$ & $7.66 \times 10^{-4}$ & $8.20 \times 10^{-4}$ & $1.87 \times 10^{-3}$ & $1.88 \times 10^{-3}$ \\
		768 & $1.46 \times 10^{-5}$ & $5.80 \times 10^{-5}$ & $6.08 \times 10^{-5}$ & $1.72 \times 10^{-4}$ & $1.73 \times 10^{-4}$ \\
		3072 & $1.00 \times 10^{-6}$ & $3.74 \times 10^{-6}$ & $3.91 \times 10^{-6}$ & $1.27 \times 10^{-5}$ & $1.27 \times 10^{-5}$ \\
		6144 & $2.59 \times 10^{-7}$ & $9.41 \times 10^{-7}$ & $9.82 \times 10^{-7}$ & $3.37 \times 10^{-6}$ & $3.37 \times 10^{-6}$ \\
		12288 & $6.68 \times 10^{-8}$ & $2.36 \times 10^{-7}$ & $2.46 \times 10^{-7}$ & $8.91 \times 10^{-7}$ & $8.91 \times 10^{-7}$ \\
		\midrule
		& $O({\alpha_1})$ & $O({\rho_1})$ & $O({\rho_2})$ & $O({u_1})$ & $O({u_2})$ \\
		\midrule
		& $1.86$ & $1.86$ & $1.88$ & $1.72$ & $1.72$ \\
		& $1.93$ & $1.98$ & $1.98$ & $1.88$ & $1.88$ \\
		& $1.95$ & $1.99$ & $1.99$ & $1.91$ & $1.91$ \\
		& $1.95$ & $2.00$ & $2.00$ & $1.92$ & $1.92$ \\
		\bottomrule
	\end{tabular}
	\label{tab:conv}
	\end{center}
\end{table}	
\subsection{Two-dimensional circular explosion problem}
To carry out another two-dimensional test, a circular explosion problem is considered. This problem can be interpreted as the two-dimensional extension of Riemann problems in radial symmetry, see also \cite{Toro2009RS}. The computational domain contains in the center a circle of radius $R=0.5$ that divides the domain into two different states, the inner state, and the outer state, defining the initial condition as follows,
\begin{equation}
	\bQ(\bx,t) = \left\{
	\begin{array}{ll}
		\bQ_L & \mbox{if }|\bx|<0.5,\\
		\bQ_R & \mbox{otherwise },
	\end{array}
	\right.
	\label{eq:CE}
\end{equation}
where $\bQ_L$ and $\bQ_R$ are described in Table~\ref{tab:CE}.
\begin{table}[!ht]
	\centering
	\caption{Left and right states of the circular explosion problem given by Equation~\eqref{eq:CE}}
	\begin{tabular}{cccccccccc}
		\toprule
		& $\alphaI$ & $\rhoI$ & $\rhoII$ & $\uI{1}$ & $\uI{2}$ & $\uI{3}$ & $\uII{1}$ & $\uII{2}$ & $\uII{3}$ \\	\midrule 
		$\bQ_L$ & 0.4 & 2 & 1.5 & 0 & 0 & 0 & 0 & 0 & 0\\
		$\bQ_R$ & 0.8 & 1 & 0.5 & 0 & 0 & 0 & 0 & 0 & 0 \\
		\bottomrule
	\end{tabular}
	\label{tab:CE}
\end{table} 
As proposed in~\cite{Dumbser2010Force2}, a reference solution can be derived by solving an equivalent PDE in radial direction with geometric source terms. To obtain more details about the radial system, the reader is referred to~\cite{Rio-Martin2023}. The numerical solution has been computed in the domain $\Omega=[-1,1]\times[-1,1]$, using the curl-free approach that has been described in this paper, considering a Cartesian staggered mesh discretized with $4800\times4800$ cells. As reported in~\cite{Rio-Martin2023}, the reference solution has been computed using a second-order TVD finite volume method with the Rusanov flux on a mesh with 128000 cells in 1D for the equivalent system in radial direction. The final time of the simulation is $t=0.1$ and the EOS for both phases is~\eqref{eq:energy_ideal_gas}, with parameters $\gammaI = 1.4$, and $\gammaII = 2$.

\begin{figure}[ht!]
	\begin{tabular}{lll} 
		\includegraphics[trim=10 15 60 50,clip,width=0.42\textwidth]{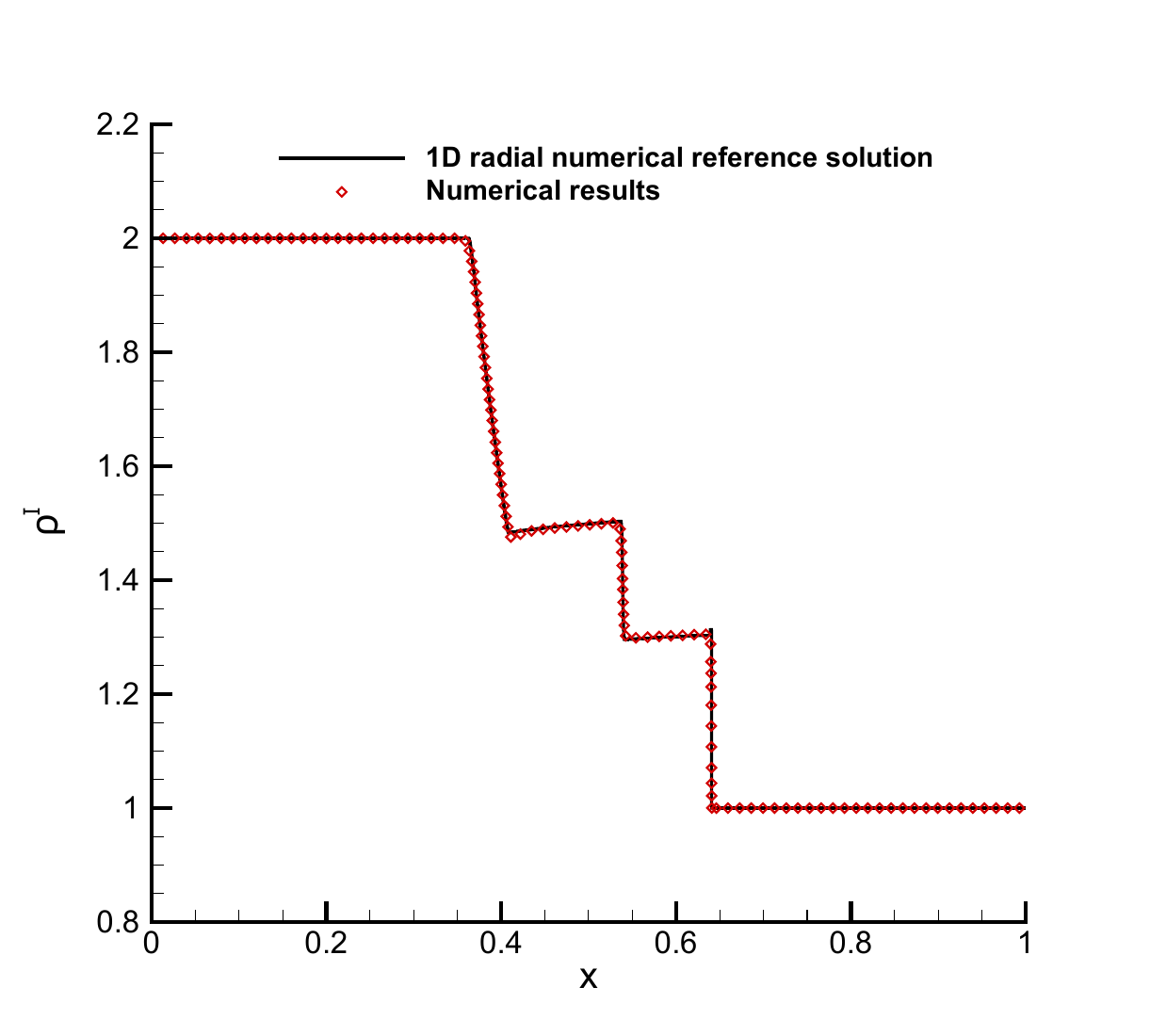}   & 
		\includegraphics[trim=10 15 60 50,clip,width=0.42\textwidth]{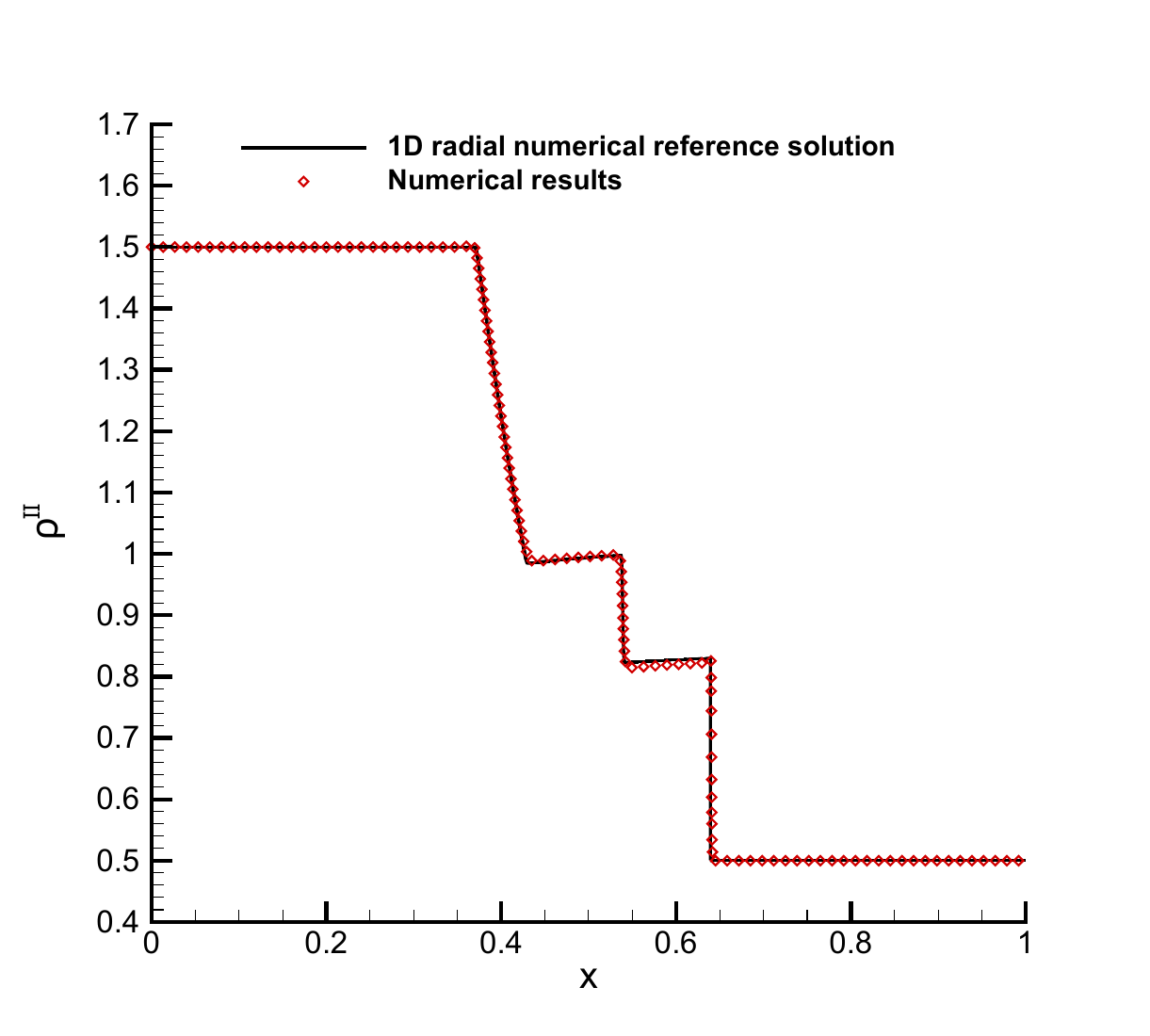}   \\ 
		\includegraphics[trim=10 15 60 50,clip,width=0.42\textwidth]{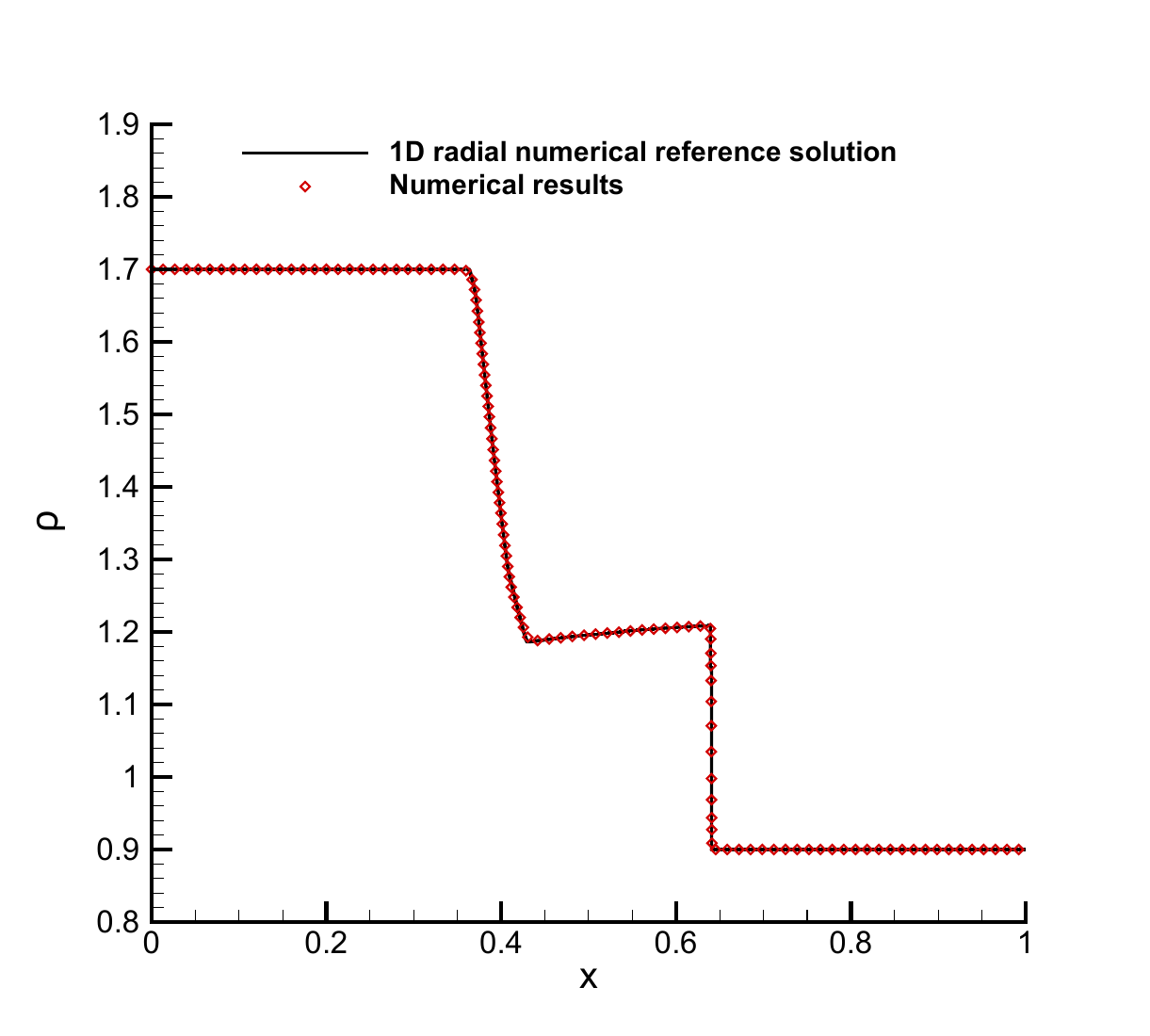}   &   
		\includegraphics[trim=10 15 60 50,clip,width=0.42\textwidth]{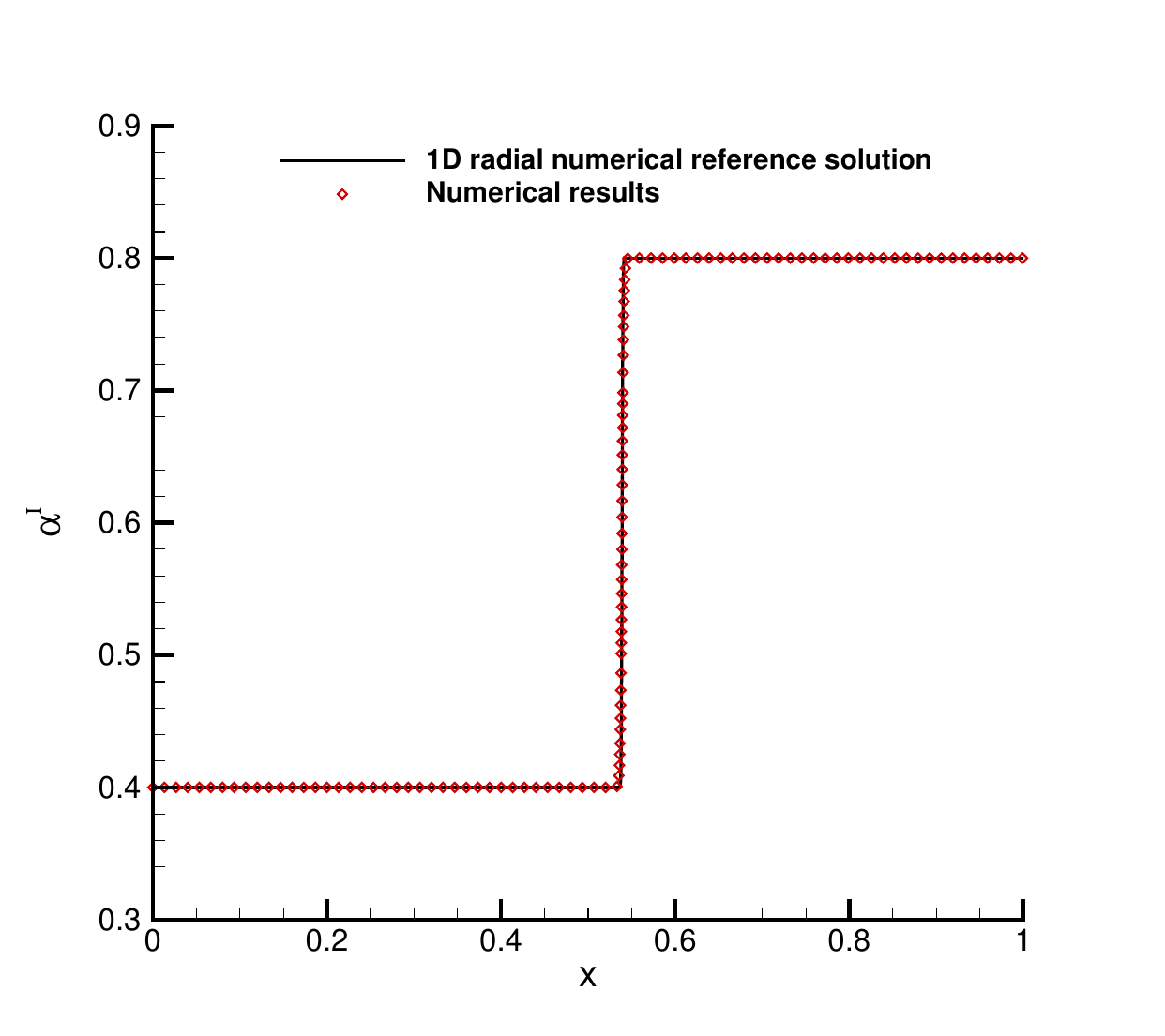}  \\
		\includegraphics[trim=10 15 60 50,clip,width=0.42\textwidth]{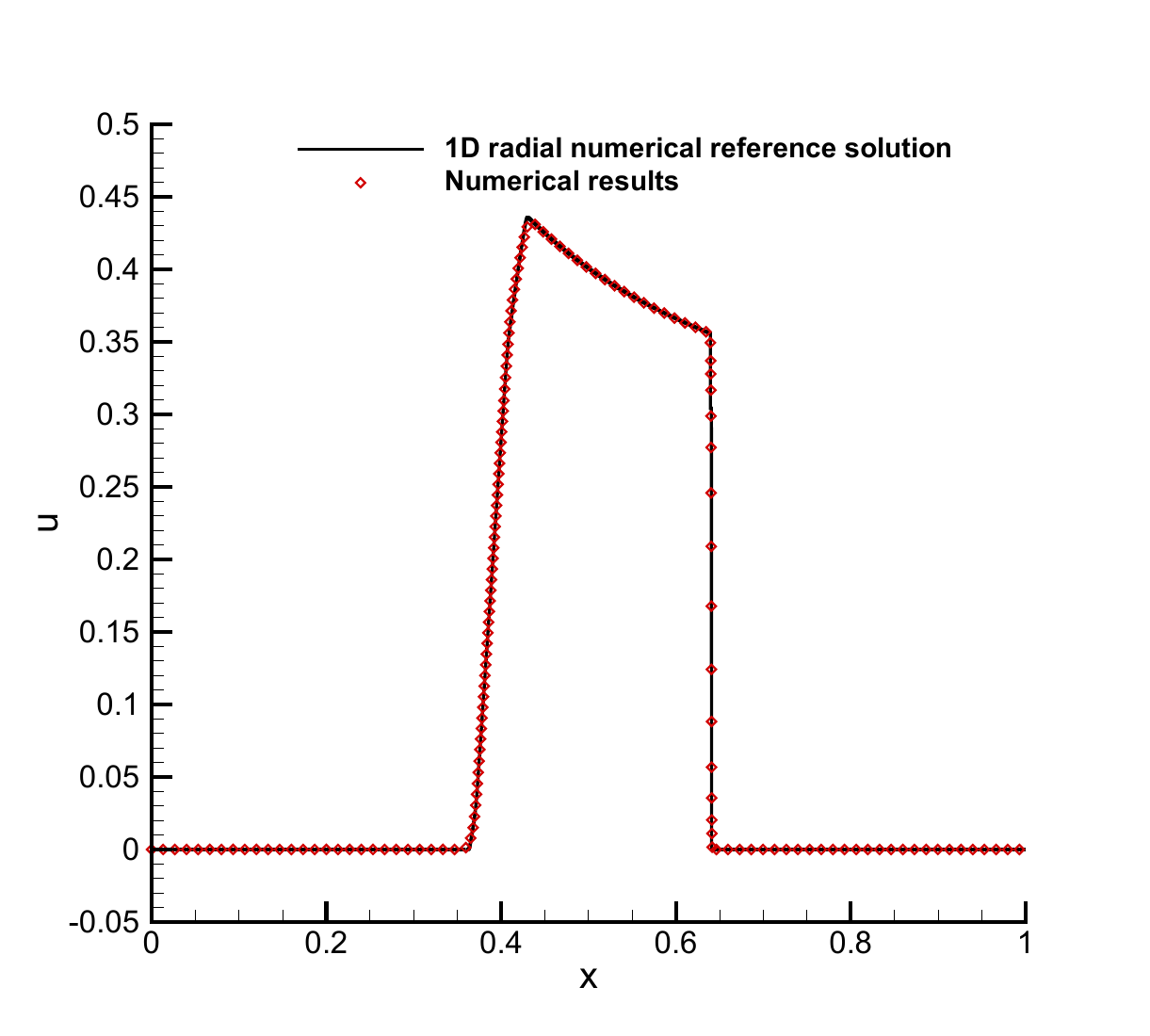}      &
		\includegraphics[trim=10 15 60 50,clip,width=0.42\textwidth]{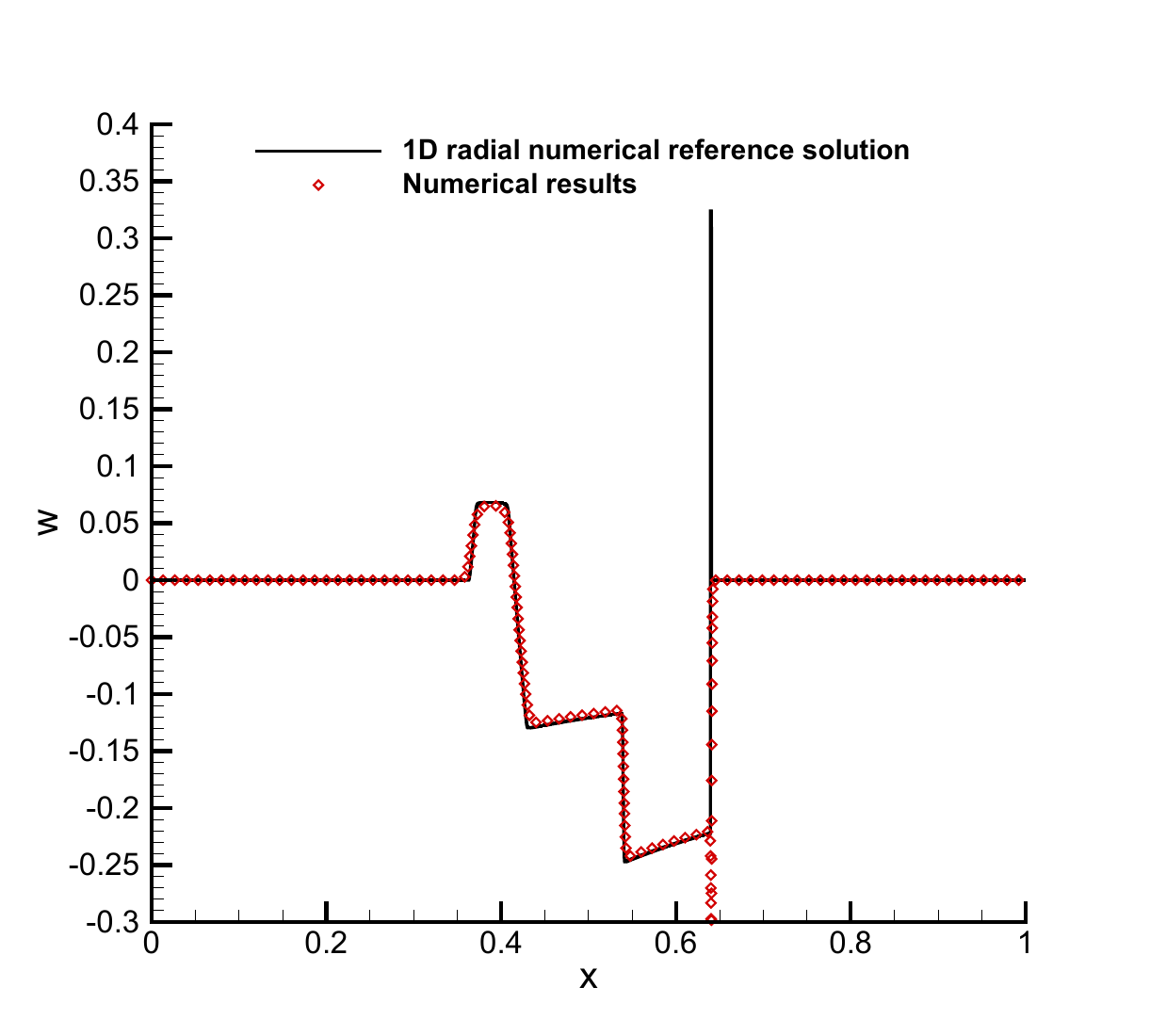}      
	\end{tabular} 
	\caption{Solution of the 2D circular explosion problem with the structured preserving finite volume scheme for the initial condition showed in Table~\ref{tab:CE} solved on a Cartesian staggered mesh with $4800\times4800$ cells at time $t=0.1$, in comparison with the radial reference solution. Top: densities of each phase, $\rhoI$ and $\rhoII$. Center: mixture density $\rho$ and volume fraction $\alphaI$. Bottom: mixture velocity $\u{}$ and relative velocity $\w{}=\uI{}-\uII{}$.} 
	\label{fig:CE}
\end{figure}

\begin{figure}[ht!]
	\begin{center}
		\includegraphics[trim=1 1 1 1,clip,width=0.65\textwidth]{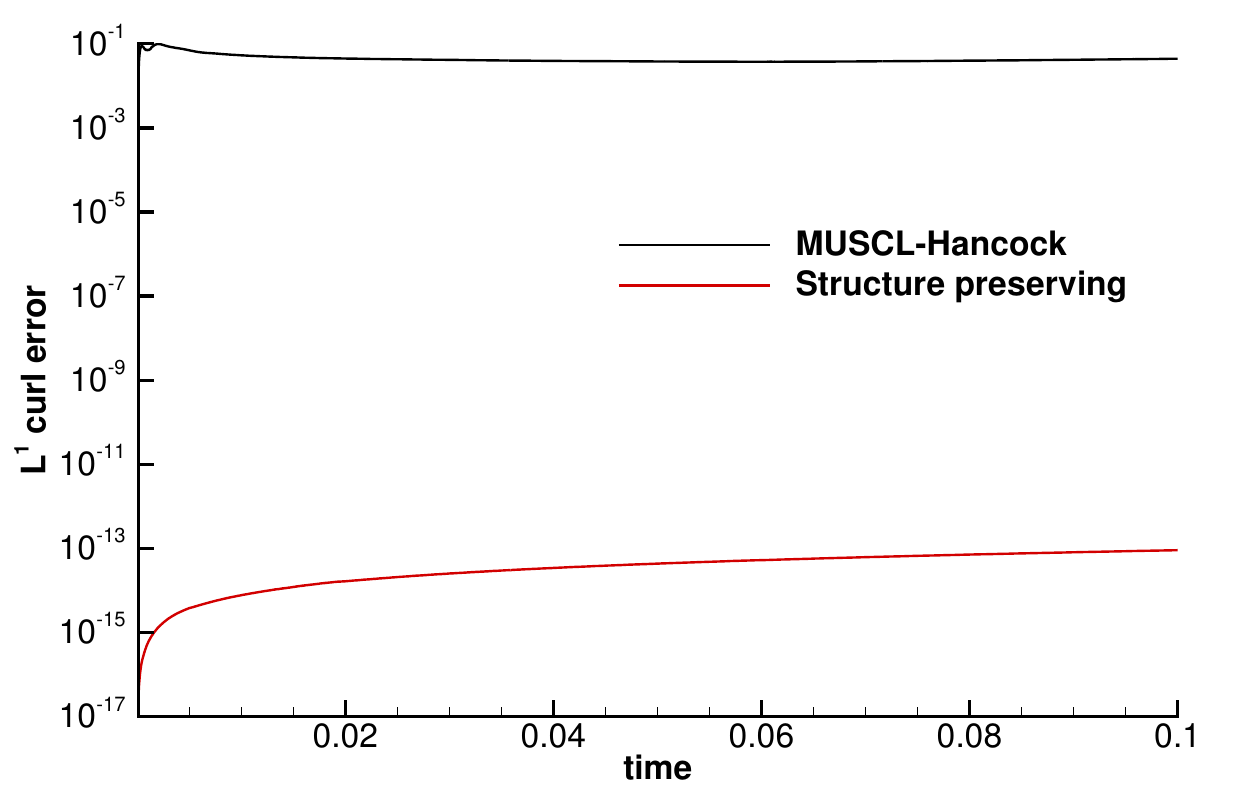}   
		\caption{Time-evolution of the $L^1$ norm of the discrete curl errors using the staggered compatible curl-free discretization (red line) and without the compatible curl-free discretization (black line) for the 2D circular explosion problem, using a mesh with $4800\times 4800$ cells.} 
		\label{fig:curlErrors_CE}
	\end{center}
\end{figure}

Figure~\ref{fig:CE} shows the numerical results compared with the reference solution described before. Moreover, Figure~\ref{fig:curlErrors_CE} illustrates the time evolution of the $L^1$ norm of the curl errors for the solution computed without the curl-free approach and for the solution calculated with the structure-preserving method detailed in this paper.

\subsection{Dambreak}
\begin{figure}[h!]
	\begin{center}
		\includegraphics{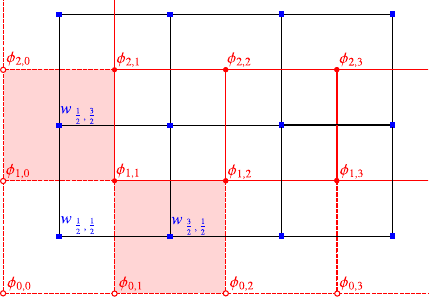}
		\caption{Schematic showing the extra layer of points of the scalar field from which $\bw$ is calculated at the boundaries. The highlighted areas are examples of stencils needed to compute $\bw$ on the $x=0$ and $y=0$ walls.}
		\label{fig:BC_dambreak}
	\end{center}
\end{figure}

\begin{figure}[h!]
	\begin{center}
		\includegraphics[trim=2 0.2 2 2,clip,width=0.5\textwidth]{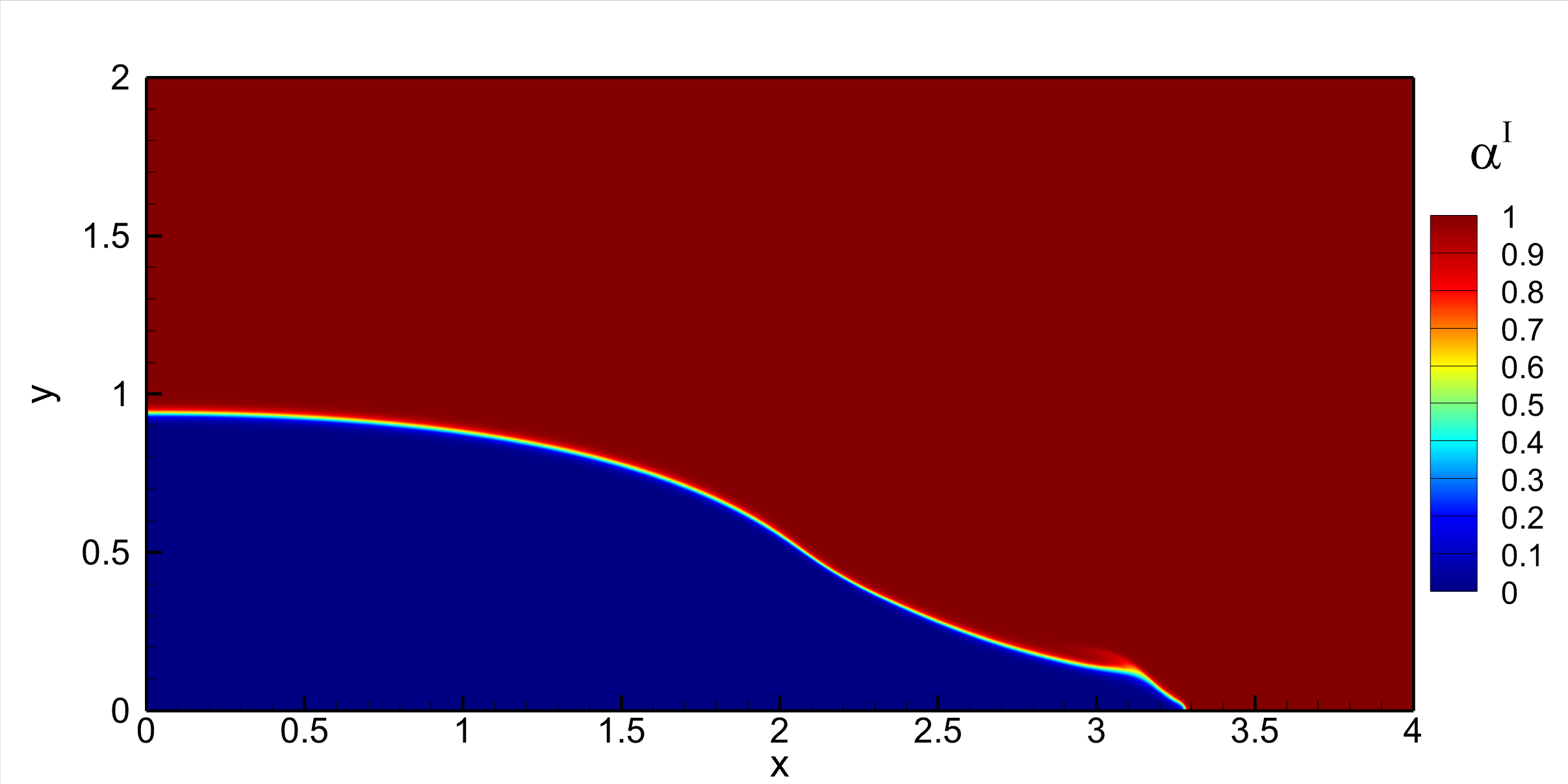}   
		\includegraphics[trim=2 0.2 2 2,clip,width=0.47\textwidth]{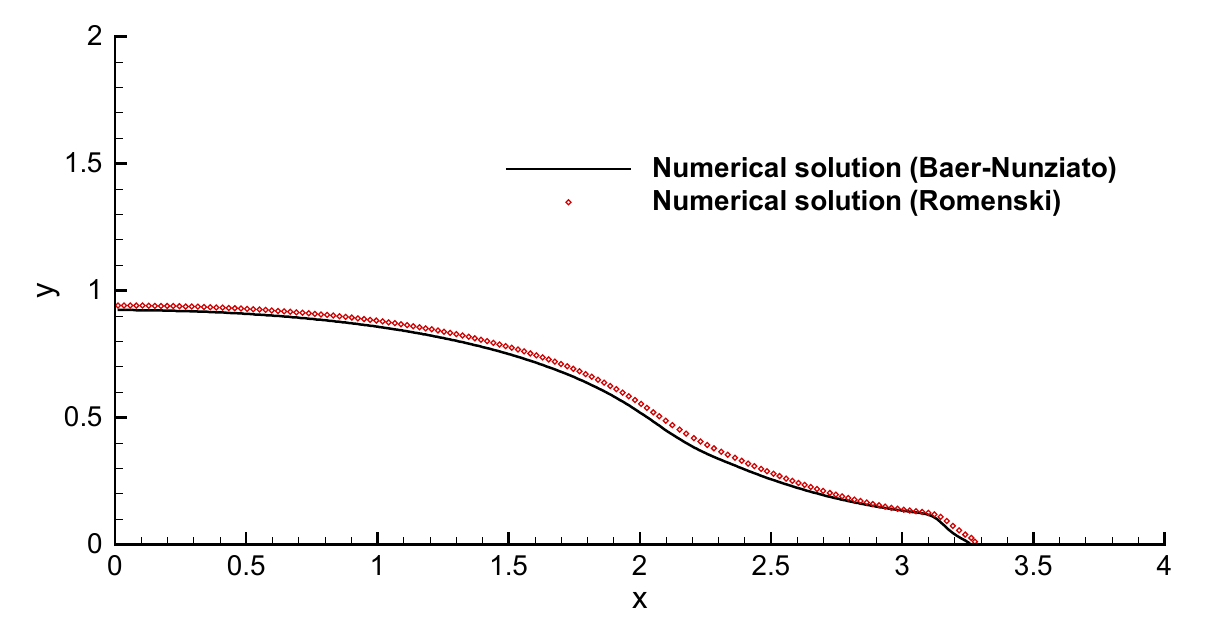}
		\caption{Dambreak problem at time $t = 0.4$. Left: Numerical solution, computed using the structure-preserving finite volume method to solve the barotropic SHTC model proposed in this paper, considering a mesh with $4800\times 2400$ cells. Right: comparison between the same numerical solution  with that of the barotropic reduced Baer--Nunziato model presented in ~\cite{Dumbser2011TwoPhase} }
		\label{fig:dambreak}
	\end{center}
\end{figure}

We consider here a 2D Riemann problem consisting in a dambreak test case. The total computational domain is $\Omega=[0,4]\times[0,2]$ discretized over $4800\times2400$ computational elements. The initial condition is such that water (assumed here as fluid $I\!I$) occupies a rectangular region in the bottom left of the domain denoted by $\Omega_W =[0,2]\times[0,1]$ while the remaining region is occupied by air. The model variables are then initialized as follows 
\begin{align*}
	(x,y) \in \Omega_W:	\ 
	\begin{cases}
		\alphaI(x,y) = 1, \\
		\rhoI(x,y)   = \prn{\rhoI_0 \, g\, (y-2)}^{1/\gammaI}\!\!\!\!\!\!, \\
		\rhoII(x,y)  = \rhoII_0, \\
		\buI{} = \buI{} = (0,0,0)^T\!\!\!\!, 
	\end{cases} 
	\!\!\text{else}:
	\begin{cases}
		\alphaI(x,y) = 0, \\
		\rhoI(x,y)   = \rhoI, \\
		\rhoII(x,y)  = \rhoII_0\prn{1+\rhoII_0 \, g\, (y-1)}^{1/\gammaII}\!\!\!\!\!\!, \\
		\buI{} = \buI{} = (0,0,0)^T\!\!\!\!. 
	\end{cases}
\end{align*}
The initial profiles of densities are such that initial pressure is hydrostatic in both phases, and the system is at rest. Water then flows under the effects of gravity, assumed to act along the $y-$axis such that $\bg=(0,-g)$ and $g=9.80$. As for the equations of state, we assume that water behaves as a stiffened gas (following the EOS~\eqref{eq:energy_stiffened}, with $\pII_0=1, \ \cII_0=20, \ \rhoII_0=1000, \gammaII=2$ ) while the surrounding air is an ideal gas ($\gammaI = 1.4$). The CFL number is fixed at $0.5$, and the final simulation time is $t=0.4$. 
Last, it is important here to discuss how the boundary conditions are implemented. Indeed, we would like to impose slip wall boundary conditions on all the domain boundaries. In terms of the mixture velocity $\bu$ and the relative velocity $\bw$, this amounts to imposing the components normal to the boundary to zero. While for the former, this can be done trivially, the uncareful treatment of this boundary condition for $\bw$ may result in the loss of the curl-free property. One solution to this problem would be introducing an additional layer of points for the scalar field from which $\bw$ is constructed, thus extending it beyond the boundaries, see Figure \ref{fig:BC_dambreak}. We can then   
impose a boundary condition directly on the scalar field, in such a way that the required boundary condition on $\bw$ is recovered, when computed from the compatible gradient.
Thus, for no-slip wall boundary conditions, we impose for the scalar field 
\vspace{0.25cm}
\begin{equation*}
	\begin{cases}
		\phi_{0,q} = \phi_{1,q}, \\
		\phi_{N_{x}+1,q} = \phi_{N_x,q}, 
	\end{cases} \forall q\in\{0..N_y\}, \quad 
	\begin{cases}
		\phi_{p,0} = \phi_{p,1}, \\
		\phi_{p,N_{y}+1} = \phi_{p,N_y},
	\end{cases}
	\forall p\in\{0..N_x\}.
\vspace{0.25cm}
\end{equation*}
The boundary conditions for the rest of variables are treated as for a classical finite volume scheme (slopes are set to zero for all variables except for the mixture velocity component normal to the boundary which is taken with opposite sign). The numerical results for this test case are presented in Figure \ref{fig:dambreak}. The left side of the figure shows the 2D distribution of the volume fraction $\alphaI$, where blue corresponds to water and red to air. The right side of the figure shows a comparison of the water depth profile extracted from this simulation, with a reference solution computed with a third-order ADER-WENO finite volume scheme on a very fine uniform Cartesian grid, solving the inviscid and barotropic reduced Baer--Nunziato model presented in~\cite{Dumbser2011TwoPhase}. The comparison shows an excellent agreement between both solutions, computed from two different models using two different numerical methods. 

\begin{figure}[!h]
	\begin{center}
		\begin{tabular}{cc} 
			\includegraphics[trim=5 10 15 10,clip,width=0.6\textwidth]{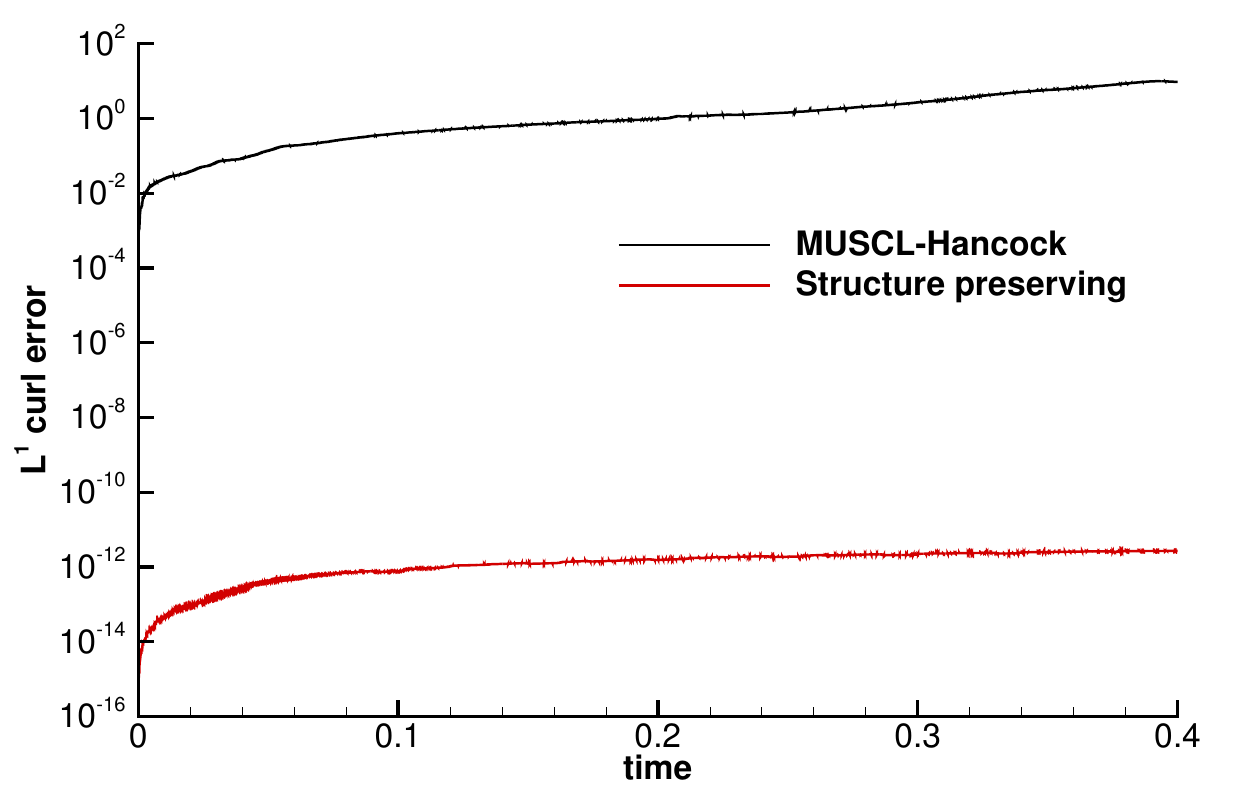}   
		\end{tabular} 
		\caption{Time-evolution of the $L^1$ norm of the discrete curl errors using the staggered compatible curl-free discretization and without the compatible curl-free discretization for the dambreak problem, using a mesh with $4800\times 2400$ cells.} 
		\label{fig:curlErrors_dambreak}
	\end{center}
\end{figure}

\subsection{Kelvin-Helmholtz instability}

\begin{figure}[h!]
	\begin{tabular}{lll}
		\includegraphics[trim=100 75 100 150,clip,width=0.47\textwidth]{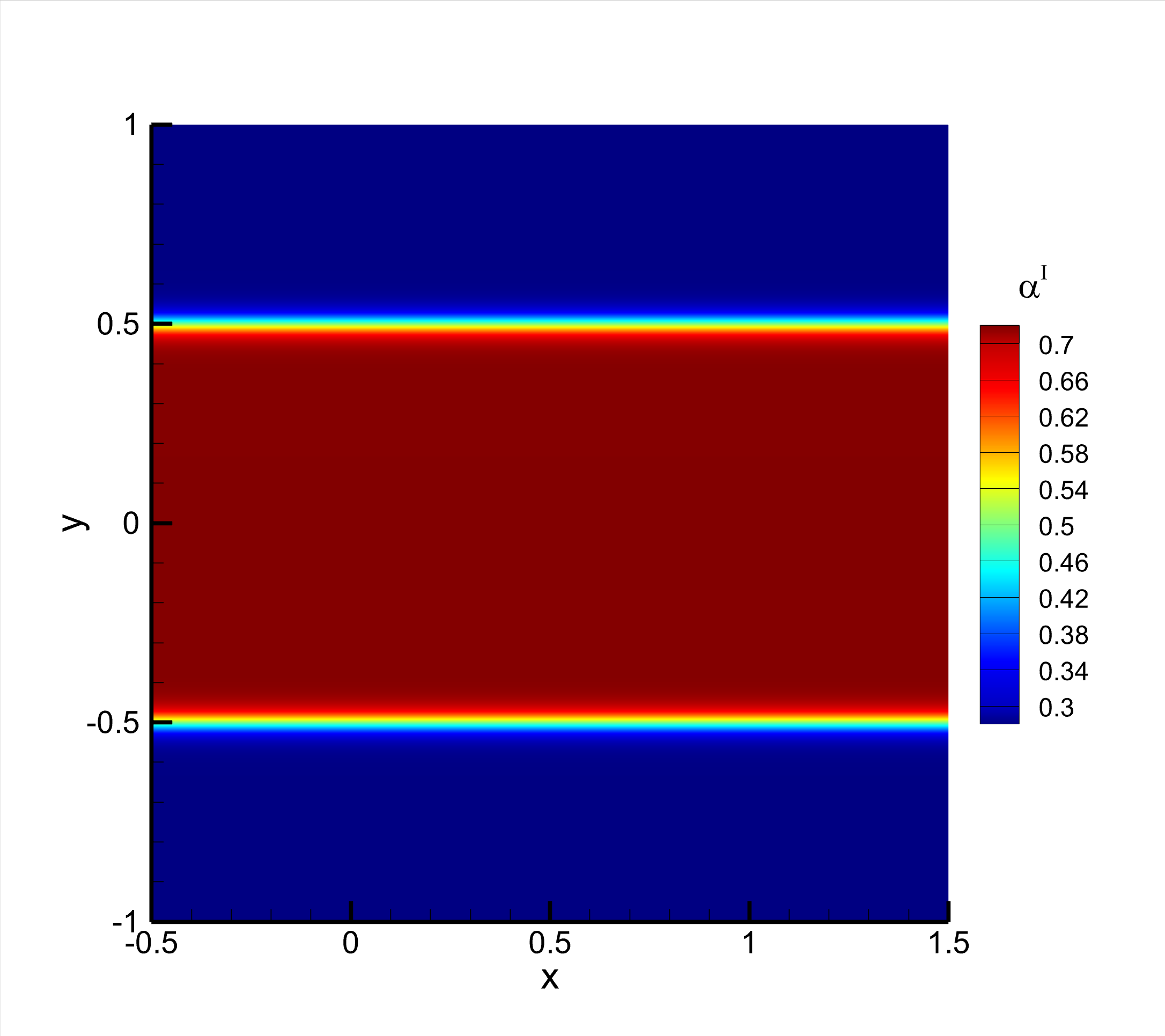}   &  
		\includegraphics[trim=100 75 100 150,clip,width=0.47\textwidth]{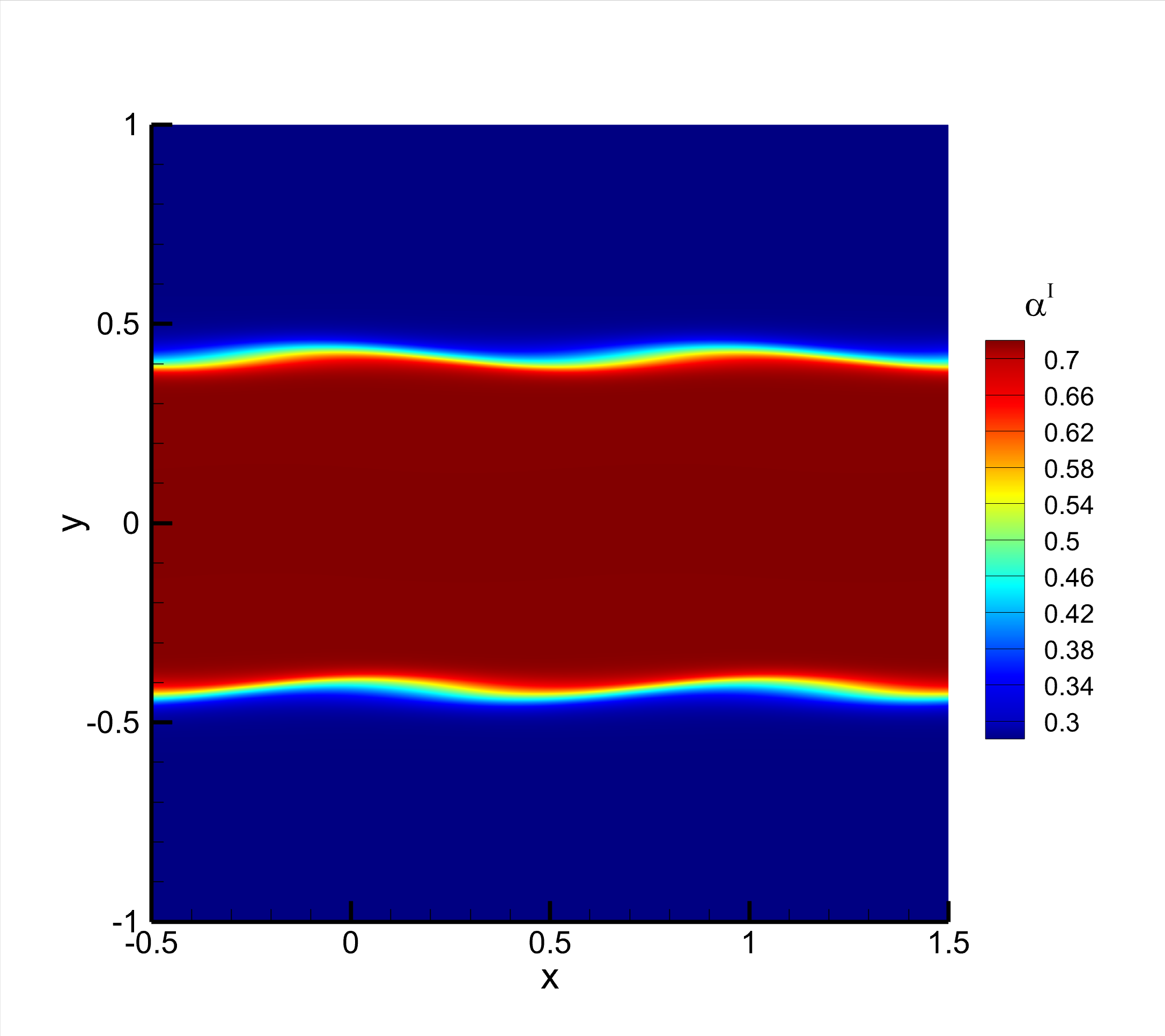}   \\
		\includegraphics[trim=100 75 100 150,clip,width=0.47\textwidth]{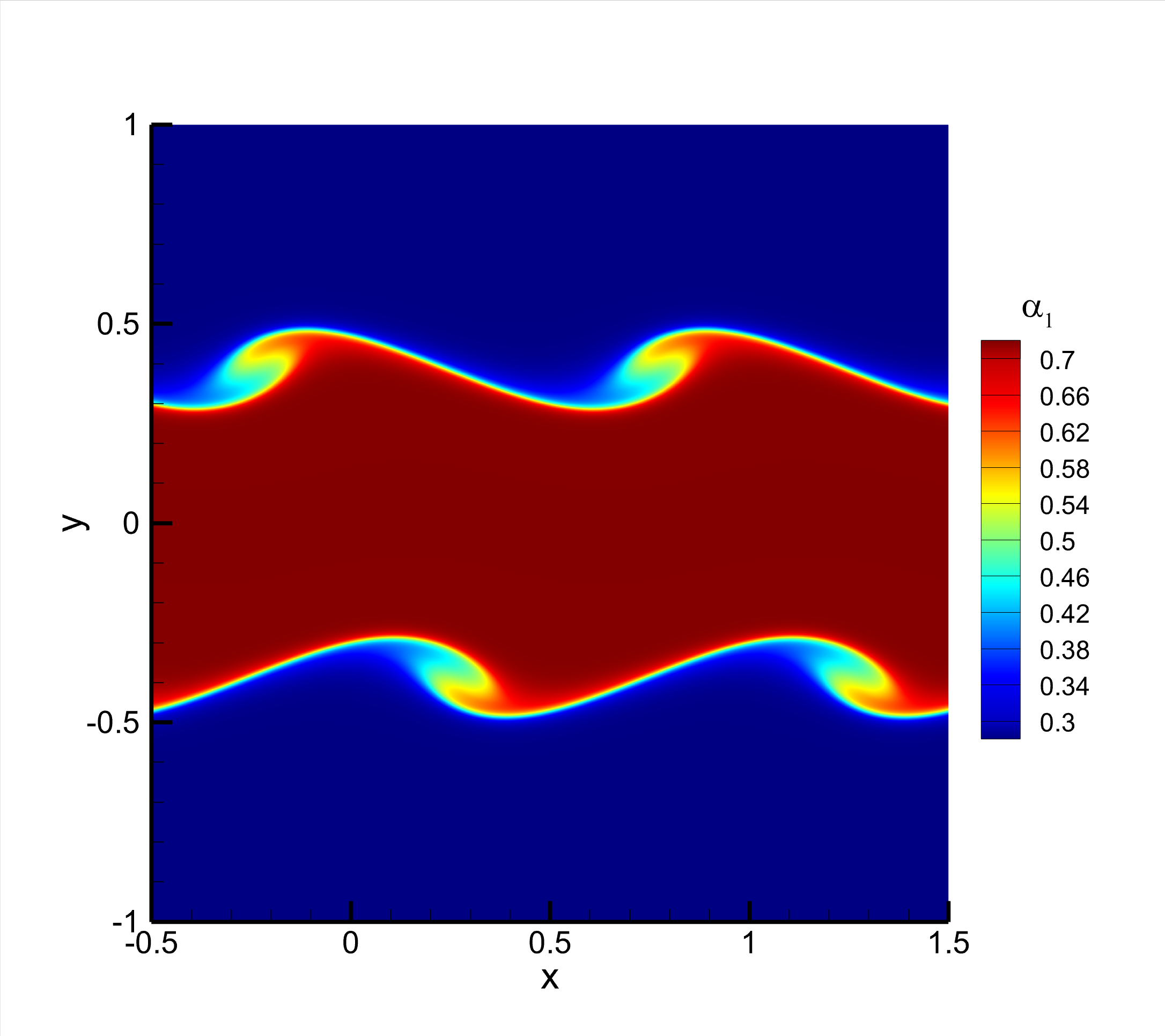}   &  
		\includegraphics[trim=100 75 100 150,clip,width=0.47\textwidth]{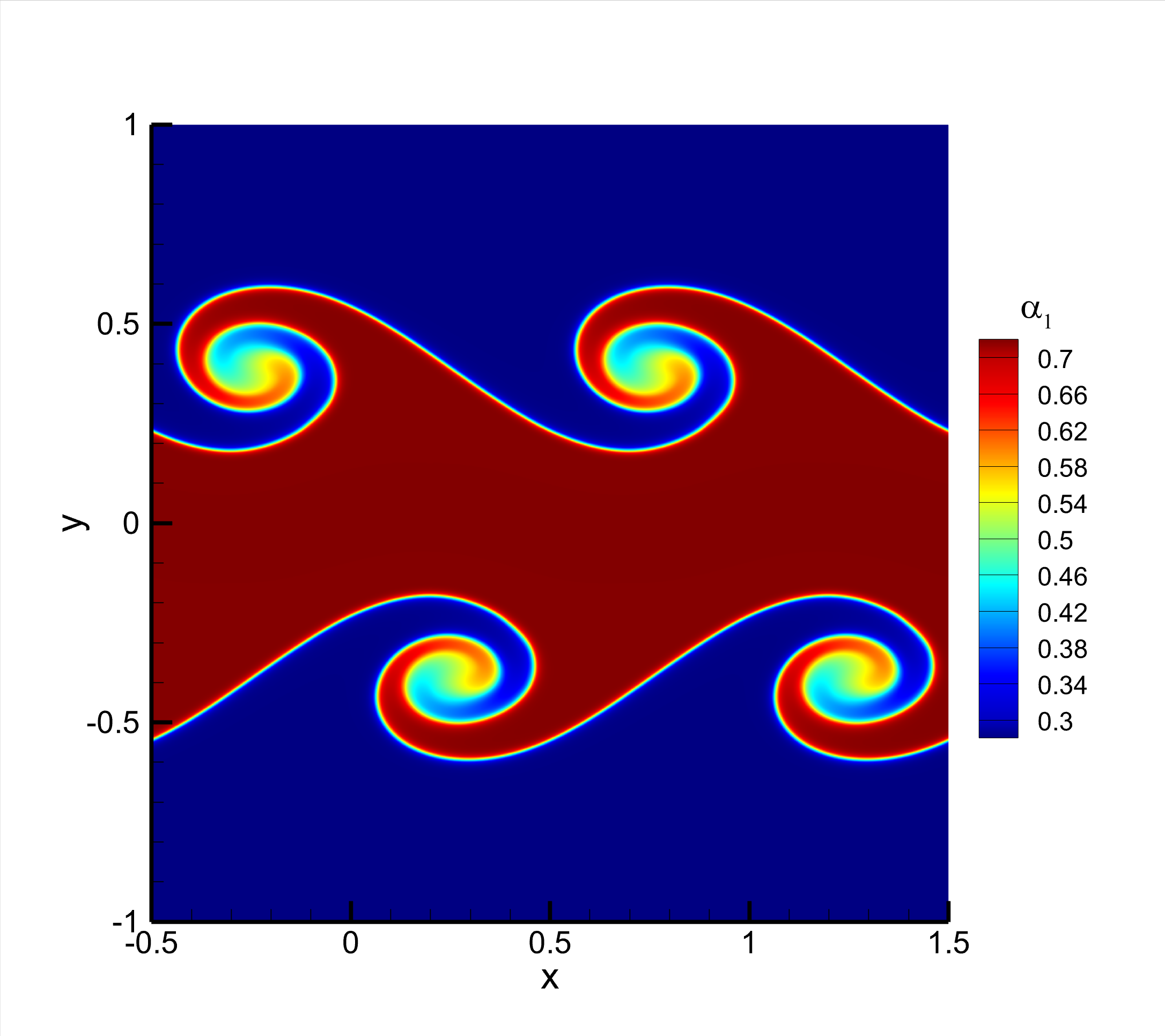}   \\   
		\includegraphics[trim=100 75 100 150,clip,width=0.47\textwidth]{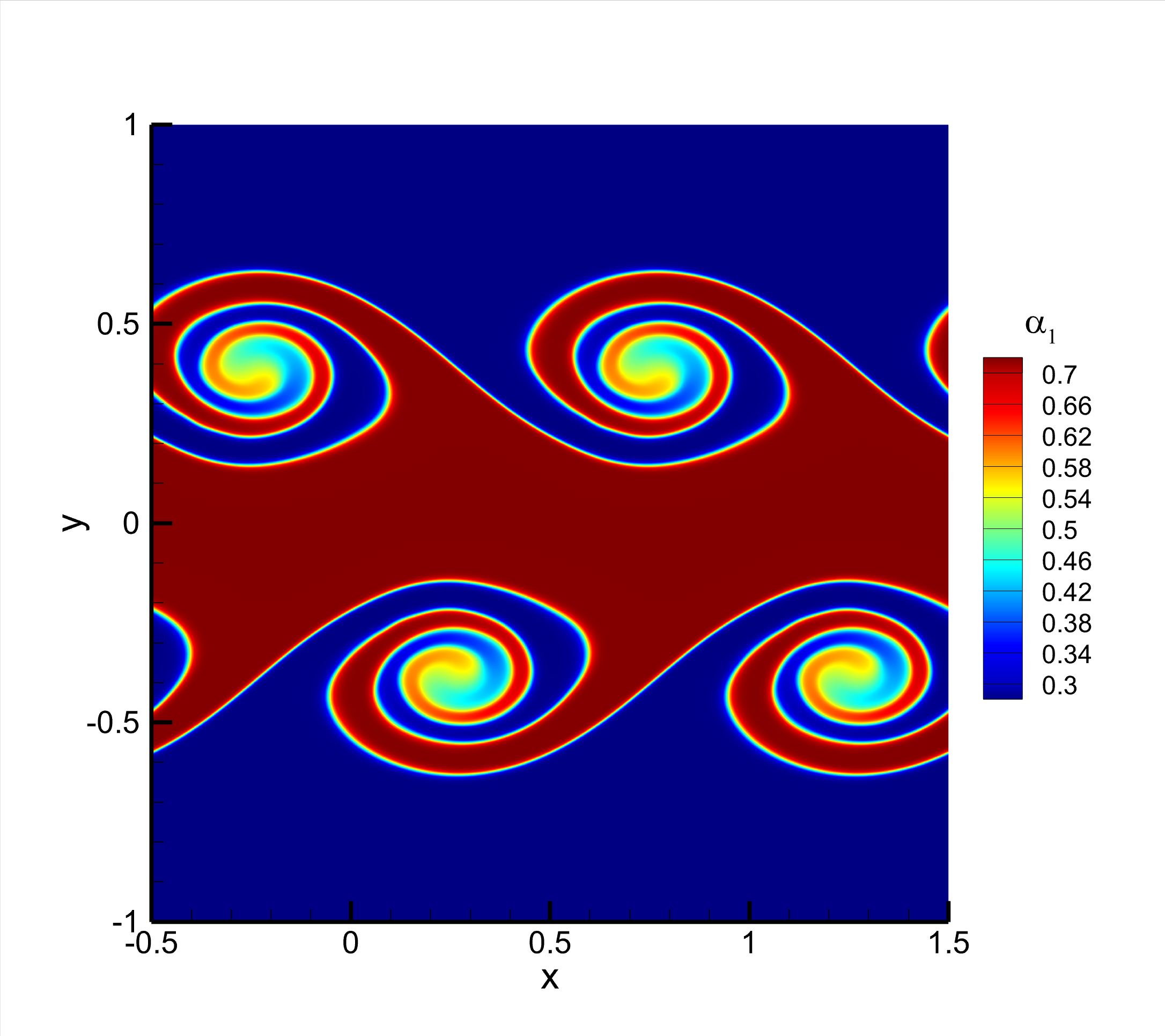}   &
		\includegraphics[trim=100 75 100 150,clip,width=0.47\textwidth]{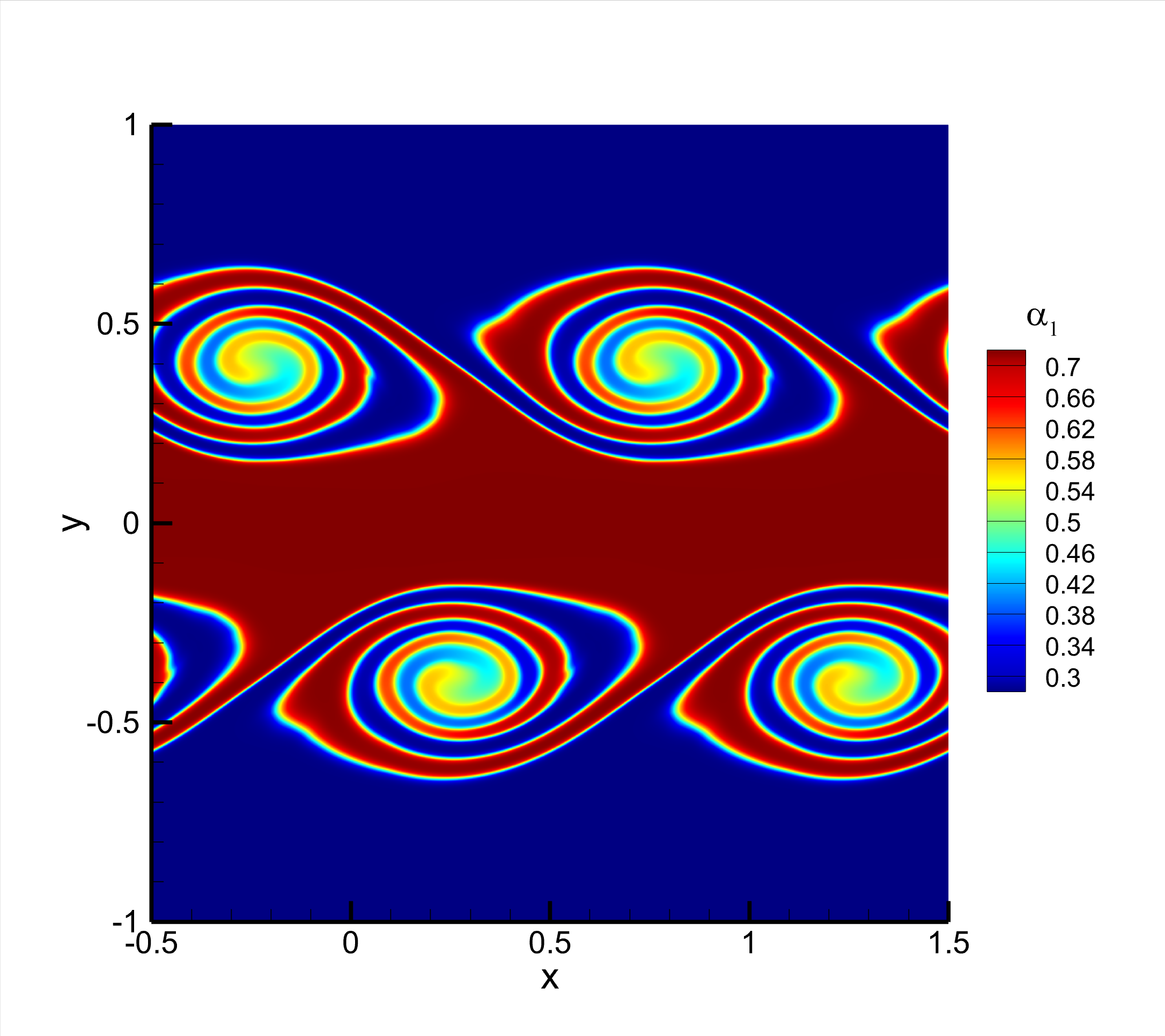}
	\end{tabular}
	\caption{Two-phase Kelvin-Helmholtz instability: numerical solution obtained with the structure-preserving finite volume scheme at times $t=0, 2, 3, 4, 5, 6$ using a mesh of $1024\times2048$ for the computational $\Omega=[-0.5,0.5]\times[-1,1]$. The plot shown in $[0.5,1.5]\times[-1,1]$ is the juxtaposition of the original numerical solution.} 
	\label{fig:KH}
\end{figure}

\begin{figure}[!h]
	\begin{center}
		\begin{tabular}{cc} 
			\includegraphics[trim=5 10 15 10,clip,width=0.6\textwidth]{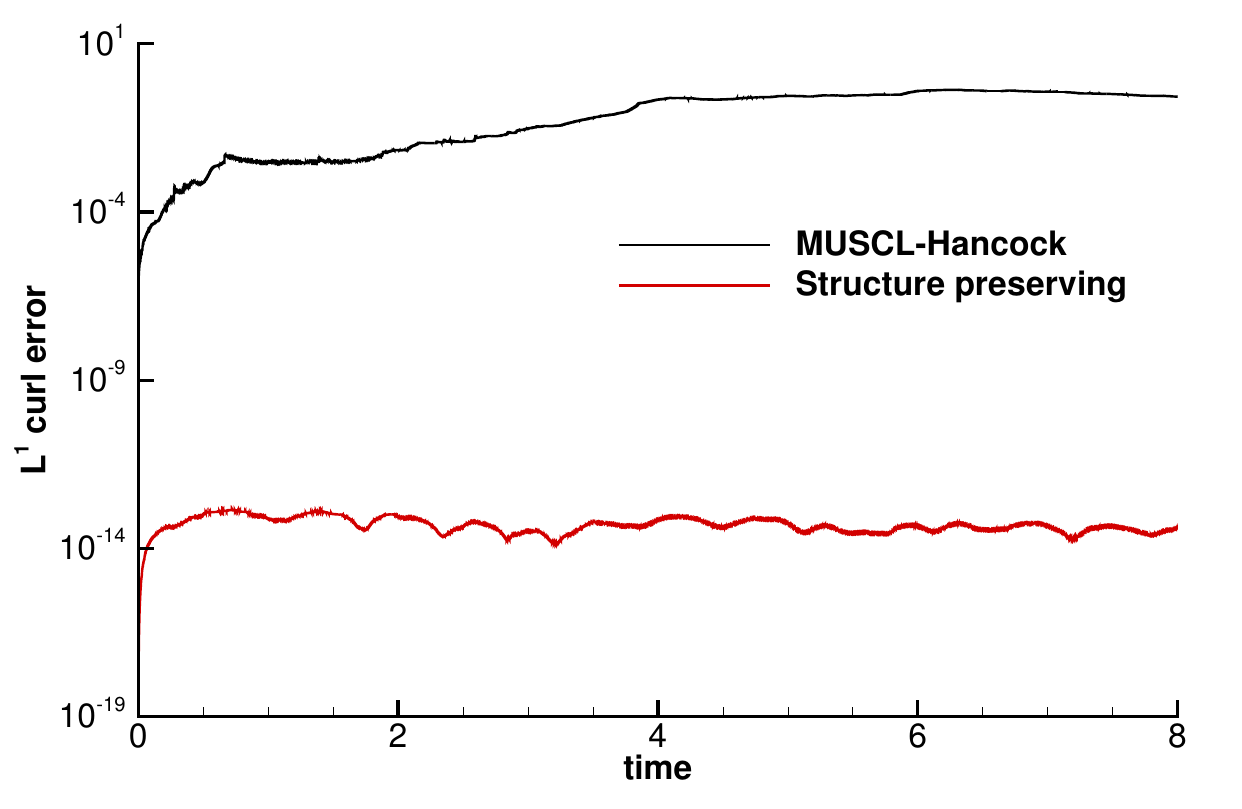}   
		\end{tabular} 
		\caption{Time-evolution of the $L^1$ norm of the discrete curl errors using the staggered compatible curl-free discretization and without the compatible curl-free discretization for the Kelvin-Helmholtz instability, using a mesh with $2048\times 2048$ cells.} 
		\label{fig:curlErrors_KH}
	\end{center}
\end{figure}

Kelvin-Helmholtz instabilities are turbulent flow patterns formed at the interface between two moving flows with different velocities or densities, which appear when a velocity gradient is created at the interface between the two fluids. In this last test case, we will perform the simulation considering two flows with different initial velocities and densities in the computational domain $\Omega=[-0.5,0.5]\times[-1,1]$. We will impose periodic boundary conditions in both $x-$ and $y-$direction and introduce a perturbation in the system to trigger instabilities at the interface between the two fluids. The simulation will be performed on a $1024\times2048$ mesh, up to a final time of $t=6$.
We will consider that at the initial time, the densities are constant $\rhoI=1$ and $\rhoII=2$. The initial condition for $\alphaI$ is given by
\begin{equation*}
	\alphaI = \left\{
	\begin{array}{ll}
		0.5+0.25\,\tanh(25(y+0.5))) & \quad\mbox{  if  } y < 0,	\\
		0.5-0.25\,\tanh(25(y-0.5))) & \quad\mbox{  if  } y\ge 0.
	\end{array}\right.
\end{equation*}
The EOS chosen for both phases is ideal gas~\eqref{eq:energy_ideal_gas}, with $\gammaI=1.4$ and $\gammaII=2$. The pressure in both phases is initialized as $\pI=\pII=100/\gammaI$, and the initial velocities are given by 
\begin{equation*}
	\uI{1}=\uII{1} = \left\{
	\begin{array}{ll}
		0.5\,\tanh(25(y+0.5))) & \quad\mbox{  if  } y < 0,\\
		-0.5\,\tanh(25(y-0.5))) & \quad\mbox{  if  } y \ge 0,\\
	\end{array}\right.
\end{equation*}
\begin{equation*}
	\uI{2}=\uII{2} = \left\{\begin{array}{ll}
		-10^{-2}\,\sin(2\pi x)\,\sin(2\pi(y+0.5)) & \quad\mbox{  if  } y < 0,\\
		10^{-2}\,\sin(2\pi x)\,\sin(2\pi(y-0.5))& \quad\mbox{  if  } y \ge 0,
	\end{array}\right.
\end{equation*}
The perturbation introduced in the $y-$component of both velocities causes the instabilities to be generated. In Figure~\ref{fig:KH}, the values of $\alphaI$ at times $t=0,2,3,4,5,6$ are shown, presenting the usual roll-up and vortex formation that is characteristic of Kelvin-Helmholtz instabilities, see also \cite{Thomann2023}. In Figure \ref{fig:curlErrors_KH} we show the temporal evolution of the curl error of the relative velocity field, which remains at the level of machine precision for all times when using the new structure-preserving FV scheme presented in this paper and which is several orders of magnitude larger when using a classical MUSCL-Hancock-type scheme. 

\section{Conclusion}
In this paper, we have numerically solved the barotropic two-phase model of Romenski \cite{Romenski1998,Romenski2004,Romenski2007TwoPhase,Romenski2010TwoPhase} using a new structure-preserving finite volume scheme that conserves the curl-free constraint of the relative velocity exactly also at the discrete level. This property was highlighted in a set of test cases, also when compatible wall boundary conditions on the relative velocity were imposed. One could extend the presented scheme to higher-order, similarly to what was done in \cite{Boscheri2023SPDG}. Another interesting approach under consideration, consists in the development of Thermodynamically Compatible Schemes \cite{Busto2022HTCGPR,Busto2023new,Abgrall2023}, allowing to conserve both the mathematical structure of the system while also keeping its energy conserved and its entropy production admissible at the discrete level. The latter property, in particular, was addressed for this model in \cite{Thomann2023}. Last, it would be of practical interest to extend the presented model to the dissipative case where viscosity and heat transfer are both present while safeguarding the hyperbolic structure, which can be done in the framework of SHTC equations, as in \cite{Peshkov2018SHTC,Peshkov2016}.

\section*{Acknowledgements}
	
	This research was funded by the Italian Ministry of Education, University and Research (MIUR) in the frame of the Departments of Excellence Initiative 2018--2027 attributed to DICAM of the University of Trento (grant L. 232/2016) and in the frame of the PRIN 2022 project \textit{High order structure-preserving semi-implicit schemes for hyperbolic equations}. 
	FD was also funded by NextGenerationEU, Azione 247 MUR Young Researchers – SoE line. LR, FD, and MD are members of the Gruppo Nazionale Calcolo Scientifico-Istituto Nazionale di Alta Matematica (GNCS-INdAM). The authors would like to acknowledge support from the CINECA under the ISCRA initiative, for the availability of high-performance computing resources and support (project number IsB27\_NeMesiS).
	This research was also co-funded by the European Union NextGenerationEU (PNRR, Spoke 7 CN HPC). Views and opinions expressed are however those of the author(s) only and do not necessarily reflect those of the European Union or the European Research Council. Neither the European Union nor the granting authority can be held responsible for them.

\bibliographystyle{abbrv}
\bibliography{./references}

\end{document}